\documentclass{article}

\usepackage{graphicx}
\usepackage[german,english]{babel}
\usepackage{amsmath}
\usepackage{amssymb}
\usepackage{amsfonts}
\usepackage{bbm}
\usepackage{url}
\usepackage{booktabs}
\usepackage{tabularx}
\usepackage{epsfig,subfigure}
\usepackage{longtable}
\usepackage{lscape}
\usepackage{psfrag}

\newtheorem{defn}{Definition}

\newtheorem{thm}[defn]{Theorem}
\newtheorem{cor}[defn]{Corollary}
\newtheorem{conj}[defn]{Conjecture}

\makeatletter

\graphicspath{{.}}

\title{Triangulated Manifolds with Few Vertices: Vertex-Transitive Triangulations I}

\author{\Large Ekkehard G.\ K\"ohler and Frank H.~Lutz}
\date{}

\begin{document}

\selectlanguage{english}

\maketitle

\bigskip
\bigskip
\bigskip

\vfill


\noindent
The computer enumeration of triangulated surfaces and combinatorial $3$-mani\-folds
was started by Altshuler and Altshuler \& Steinberg in the early and mid seventies of the twentieth century:
They explicitely enumerated all combinatorial $3$-mani\-folds with up to $9$ vertices 
\cite{Altshuler1974}, \cite{AltshulerSteinberg1973},
\cite{AltshulerSteinberg1974}, \cite{AltshulerSteinberg1976}, 
all neighborly combinatorial $3$-manifolds with $10$ vertices
\cite{Altshuler1977}, and all neighborly
triangulated orientable surfaces with $12$ vertices \cite{AltshulerBokowskiSchuchert1996}.

A conceptually slightly different enumeration algorithm for neighborly combinatorial 
$3$-manifolds with a \emph{vertex-transitive} cyclic or dihedral group action was presented 
by K\"uhnel and Lassmann~\cite{KuehnelLassmann1985-di} in 1985. 
With an implementation of their algorithm they enumerated 
all vertex-transitive neighborly $3$-manifolds with 
a cyclic action for up to $15$ vertices and with a dihedral action
for up to $19$ vertices. 

Further enumerations results for $4$-manifolds 
were obtained by K\"uhnel and Lassmann~\cite{KuehnelLassmann1983-unique}
(on the uniqueness of K\"uhnel's $9$-vertex triangulation of the 
complex projective plane \cite{KuehnelBanchoff1983}),
by Lassmann and Sparla \cite{LassmannSparla2000}
(on centrally symmetric triangulations of $S^2\!\times\!S^2$
with $12$ vertices), and 
by Casella and K\"uhnel \cite{CasellaKuehnel2001}
(on the existence of a $16$-vertex triangulation of the K3-surface).

In general, however, there is no algorithm to enumerate combinatorial
manifolds of dimension $d\geq 6$: By fundamental work of 
Novikov (cf.~\cite{VolodinKuznetsovFomenko1974}),
there is no algorithm to recognize (combinatorial) triangulations of spheres
of dimension $d-1\geq 5$, which would be needed to determine whether
vertex-links are spheres (see
\cite{KuehnelLutz2003pre}
for a discussion).

In this 
paper,
we describe an algorithm for the enumeration of (candidates of) 
vertex-transitive combinatorial $d$-manifolds. 
With a GAP implementation, MANIFOLD\_VT~\cite{Lutz_MANIFOLD_VT}, 
of our algorithm, we determine, up to combinatorial
equivalence, all combinatorial manifolds with a
ver\-tex-tran\-si\-tive automorphism group on\, $n\leq 13$ vertices.
With the exception of actions of groups of small order,
the enumeration is extended to $14$ and $15$ vertices. 

\pagebreak

Our enumeration algorithm is, in part, based on 
the algorithm by K\"uhnel and Lassmann.
Improvements and variants of our algorithm
are used to enumerate all vertex-transitive triangulations of 
$3$-mani\-folds with $16$ and $17$ vertices and
all vertex-transitive neighborly surfaces with up to $22$ vertices  
\cite{Lutz2004bpre},
centrally symmetric triangulations with a vertex-transitive cyclic action
\cite{Lutz2004apre},
vertex-transitive combinatorial pseudo\-mani\-folds
\cite{Lutz2004dpre},
all triangulated surfaces with $9$ and $10$ vertices \cite{Lutz2005apre},
and all combinatorial $3$-manifolds with $10$ vertices \cite{BokowskiBremnerLutzMartin2003pre}.

\section{The Enumeration Algorithm}

The aim of this section it to describe a basic algorithm for the
enumeration of candidates for ver\-tex-tran\-si\-tive combinatorial (respectively simplicial) 
manifolds as well as heuristical steps to further analyze these candidates.

The procedure consists of seven steps. In Step~1, 
the input parameters have to be fixed: the number of vertices $n$, 
the dimension $d$, and the vertex-transitive group action~$n^i$ on $n$ vertices. 
In Step~2, all candidates for simplicial $d$-manifolds 
with $n$ vertices are generated that are invariant under the transitive
group action $n^i$. These candidates are tested heuristically
whether they are manifolds in Step~3, then classified up to combinatorial
equivalence in Step~4, and examined further in Steps~5--6. 
If they are combinatorial manifolds, then in Step~7
we heuristically try to determine their topological types.

\newcounter{null}
\begin{list}
{\Roman{null}}{\usecounter{null}\setlength{\leftmargin}{5mm}\setlength{\labelwidth}{0mm}\setlength{\labelsep}{5mm}}

\item[\emph{Step}\hspace{.25mm}] \emph{\hspace{-4.5mm}1: Fix the number of vertices\, $n\geq 4$ and
            the dimension\, $2\leq d\leq n-2$.}

\medskip

  Before starting with the enumeration of vertex-transitive triangulated
  $d$-mani\-folds with $n$ vertices we need to know all vertex-transitive group
  actions on the given number $n$ of vertices.
  The transitive permutation groups of small degree~$n\leq 31$ were classified
  by Miller \cite{Miller1896-12}, \cite{Miller1897-13}, Butler and McKay \cite{ButlerMcKay1983},
  Royle \cite{Royle1987}, Butler \cite{Butler1993}, Conway, Hulpke, and McKay \cite{ConwayHulpkeMcKay1998},
  and  Hulpke \cite{Hulpkepre}, \cite{Hulpke1996}; 
  see Table~\ref{tbl:transgroups} for the numbers of distinct actions
  that occur. In general, a finite group can have different group
  actions on $n$ vertices. The respective
  permutation groups then are non-isomorphic as permutation groups
  but isomorphic as finite groups.
  If $n$ is prime, then the corresponding transitive permutation
  groups are primitive. Generators for all the transitive
  permutation groups of degree $n\leq 31$ are available via
  the transitive permutation group library of the computer algebra
  package GAP \cite{GAP4}.


\begin{table}
\small\centering
\defaultaddspace=0.1em
\caption{The number of transitive group actions on $n\leq 31$ vertices.}\label{tbl:transgroups}
\begin{tabular*}{\linewidth}{c@{\hspace{3mm}}r@{\hspace{3mm}}r@{\hspace{3mm}}r@{\hspace{3mm}}r@{\hspace{3mm}}r@{\hspace{3mm}}r@{\hspace{3mm}}r@{\hspace{3mm}}r@{\hspace{3mm}}r@{\hspace{3mm}}r@{\hspace{3mm}}r@{\hspace{3mm}}r@{\hspace{3mm}}r@{\hspace{3mm}}r@{}}
 \addlinespace
 \addlinespace
 \addlinespace
 \addlinespace
\toprule
 \addlinespace
 \addlinespace
        $n$ &     4 & 5 &  6 & 7 &  8 &  9 & 10 & 11 &  12 & 13 &  14 &  15 &   16 & 17 \\ \midrule
 \addlinespace
 \addlinespace
 \addlinespace
        \#  &     5 & 5 & 16 & 7 & 50 & 34 & 45 &  8 & 301 &  9 &  63 & 104 & 1954 & 10 \\ \bottomrule
\\[1mm]
\toprule
 \addlinespace
 \addlinespace
        $n$ &    18 & 19 &   20 &  21 & 22 & 23 &    24 &  25 & 26 &   27 &   28 & 29 &   30 &  31\\ \midrule
 \addlinespace
 \addlinespace
 \addlinespace
        \#  &  983 &  8 & 1117 & 164 & 59 &  7 & 25000 & 211 & 96 & 2392 & 1854 &  8 & 5712 &  12\\ \bottomrule
\end{tabular*}
\end{table}

\medskip

  For every fixed pair of $n$ and $d$, we treat the corresponding transitive group actions 
  in decreasing group order. The two transitive group actions on $n$
  vertices of largest group order are the actions of the symmetric group $S_n$ 
  and of the alternating group $A_n$. 
  Both groups are transitive on the unordered $(d+1)$-subsets of the
  set $\{ 1,\dots,n\}$ for all $1\leq d\leq n-2$. Since the only $d$-manifold 
  on $n$ vertices, invariant under one of these two actions, 
  is the $(n-2)$-sphere~$S^{n-2}$ triangulated
  as the boundary of the $(n-1)$-simplex with $n=d+2$ vertices, 
  we can start for $n>d+2$ with the next smaller permutation group 
  from the list of group actions for the respective $n$. 

\medskip

  \emph{Let $n^i$ be the $i$-th transitive permutation group on $n$
  vertices from the GAP~library and let\, $\mbox{}^d\hspace{.3pt}n^i$ be
  its induced action on the $(d+1)$-subsets of the
  set $\{ 1,\dots,n\}$.}

\medskip

  Every time that our enumeration algorithm produces a new candidate for a
  combinatorial manifold in Step~2, we test in Step~4 whether we have found 
  a combinatorially equivalent candidate before.
  If not, then the combinatorial automorphism group of 
  our candidat is the current permutation group $n^i$.
  (Since we proceed with decreasing group order, the order of the
  automorphism group of our candidate cannot be larger: 
  All examples with a larger vertex-transitive group action have 
  been enumerated before.)

\bigskip

\item[\emph{Step}\hspace{.25mm}] \emph{\hspace{-5.5mm} 2 (Enumeration): Determine all pure $d$-dimensional simplicial
  complexes on $n$ vertices that have the pseudomanifold property 
  and that are invariant unter the vertex-transitive group action\, $n^i$.}

\medskip

  The `candidates' for vertex-transitive combinatorial manifolds with $n$ vertices 
  that we are going to build are pure $d$-dimensional 
  simplicial complexes $M$ that are invariant under the group
  action $n^i$. The collection of facets of every such $M$ 
  is a union of orbits of $(d\!+\!1)$-tuples with respect to 
  the induced action\, $\mbox{}^d\hspace{.3pt}n^i$\, of the permutation group $n^i$ on the set of $(d\!+\!1)$-sub\-sets
  of the ground set $\{ 1,\dots,n\}$.

\medskip

  In addition, we require that $M$ has the \emph{pseudomanifold property},
  that is, every $(d\!-\!1)$-dimen\-sion\-al \index{pseudomanifold property}
  face of $M$ must be contained in \emph{precisely two} $d$-dimensional facets.
  By transitivity, we say that an orbit of $(d\!-\!1)$-dimensional
  faces is \emph{included\, $t$ times} in an orbit of $d$-dimensional facets 
  if each $(d\!-\!1)$-dimensional member of the first orbit is included in
  $t$ sets of the latter orbit. If there is a $(d\!-\!1)$-dimensional orbit
  that is included three or more times in a $d$-dimensional orbit, then this
  $d$-dimensional orbit cannot be used for composing $M$, since this would
  violate the pseudomanifold property. 
  In a preprocessing step we sort out all
  these $d$-orbits. It then can happen that there are some $(d\!-\!1)$-orbits
  that are not included (or included only once) in any (in one) of
  the remaining $d$-orbits. We sort out these $(d\!-\!1)$-orbits 
  (and the $d$-orbits containing these $(d\!-\!1)$-orbits) as well
  and iterate this procedure as long as possible.

\medskip

  We next associate a weighted bipartite graph with the remaining 
  $d$- and $(d-1)$-orbits as nodes and an edge of weight $t$ 
  between two nodes whenever a $(d\!-\!1)$-orbit is included $t$-times 
  in a $d$-orbit. 
  Let us, for example, consider the action $7^2$ of the dihedral
  group $D_7$ on $n=7$ vertices
  and let $d=3$. There are four $3$-di\-men\-si\-on\-al orbits 

    \addtolength{\arraycolsep}{.85mm} 
  $$\begin{array}{l@{\hspace{3mm}}l}
  a:   &  1234, 2345, 3456, 4567, 1567, 1237, 1267 \\
  b:   &  1235, 2346, 2456, 3457, 1345, 3567, 1456, \\ 
       &  2347, 1467, 1247, 2567, 1236, 1257, 1367 \\
  c:   &  1245, 2356, 3467, 1457, 1347, 1256, 2367 \\
  d:   &  1246, 2357, 1356, 1346, 2457, 2467, 1357 
  \end{array}$$
  \addtolength{\arraycolsep}{-.85mm} 

  (where $1234$ denotes the tetrahedron with vertices $1$, $2$, $3$,
  and $4$, etc.)
  and four $2$-di\-men\-si\-on\-al orbits

    \addtolength{\arraycolsep}{.85mm} 
  $$\begin{array}{l@{\hspace{3mm}}l}
  e:   &  123, 234, 456, 345, 567, 167, 127 \\
  f:   &  124, 235, 356, 346, 245, 467, 457, \\
       &  134, 157, 137, 156, 237, 126, 267 \\
  g:   &  125, 236, 256, 347, 145, 367, 147 \\
  h:   &  135, 246, 357, 146, 247, 257, 136 
  \end{array}$$
  \addtolength{\arraycolsep}{-.85mm} 

  with associated weighted graph 

  \medskip

  \begin{figure}
    \begin{center}
       \mbox{}\hspace{7mm}\includegraphics{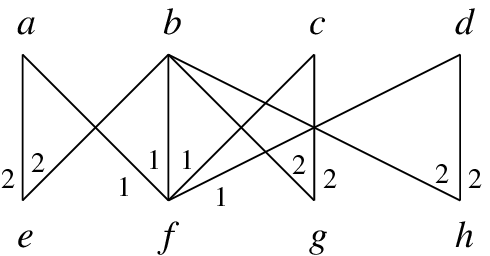}
    \end{center}
    \vspace{-2mm}
  \end{figure}

  For composing a vertex-transitive pure $d$-dimensional simplicial complex 
  with the pseudomanifold property we have to form combinations of $d$-orbits, 
  such that the resulting total weight of (the incident edges of) every 
  contained $(d\!-\!1)$-orbit is exactly two.
  To find such combinations fast, we have to choose appropriate 
  data structures for the enumeration. If the acting group $G$ has small
  group order, i.e., $|G|=m\cdot n$ with $m$ small, then the corresponding 
  weighted bipartite graph will be sparse. Thus we best use adjacency
  lists to represent the graphs. Since the graphs are bipartite,
  the resulting lists can be displayed for every graph in form of a matrix
  (with missing entries) that has a row for every $d$-orbit
  and a column for every $(d-1)$-orbit. For the above graph the matrix is 

  \addtolength{\arraycolsep}{.85mm} 
  $$\begin{array}{rl}
    \begin{array}{ccccc} e & f & g & h & \mbox{} \end{array} & \\[.5mm]
    \begin{array}{c} a\\b\\c\\d \end{array} 
    \left(
      \begin{array}{cccc}
        2 & 1 &   &   \\
        2 & 1 & 2 & 2 \\
          & 1 & 2 &   \\
          & 1 &   & 2
      \end{array}
    \right) 
    & \begin{array}{c} \\ \\ \\ \end{array}
  \end{array}$$
  \addtolength{\arraycolsep}{-.85mm} 

  In terms of this associated matrix, the problem of finding pure
  vertex-transitive simplicial complexes with the pseudomanifold
  property translates to finding all combinations of row vectors 
  such that their vector sum has entries $0$ or $2$ only.
  Missing entries `contribute' $0$ in the summation
  and therefore can be neglected in the computation.
  In our implementation we determine the valid combinations via \emph{backtracking}.

\medskip

  We group the rows of the matrix in \emph{blocks}, such that the rows of every
  block have their first non-zero entry at the same position.
  If we assume that the corresponding orbits of facets were 
  ordered lexicographically, then the rows of the first block
  have their first non-zero entry at the first position,
  and from block to block the position of the first non-zero entry
  of the respective rows increases. The above matrix has two blocks.

\medskip

  We start the backtracking with the zero row-vector as \emph{current sum vector} 
  and introduce a \emph{pointer} that points to the 
  next row-vector that is to be added to the current sum vector. 
  Initially, the pointer is set to the first row. When the first row 
  has been added to the zero vector the pointer is set to the second row as the 
  next row to be added, if possible: As soon as the current sum vector has $2$
  as its first entry, then no further row from the first block can be added  
  without violating the pseudo-manifold property.
  Thus we can set the pointer to the first row of the next block of
  the matrix, etc.

\medskip

  As soon as the current sum vector is \emph{closed}, i.e., has entries $2$ or $0$ only,
  the corresponding combination of $d$-orbits gives a vertex-transitive 
  pure simplicial complex with the pseudomanifold-property
  and therefore a candidate for a vertex-transitive simplicial manifold.
  If we would add further rows of the matrix to a closed vector,
  then we might eventually obtain another closed vector that is the sum
  of two closed vectors. However, the corresponding simplicial complex
  then is \emph{not strongly connected} (i.e., there is a pair of facets
  that cannot be joined by a path which moves from facet to facet 
  only across $(d-1)$-faces), and therefore cannot be a
  connected manifold. To avoid this, we set the pointer to END.

\medskip

  Whenever the pointer points to END, then in the following step we
  go one level up in the backtracking tree (by subtracting
  the last row of the combination from the current sum vector) 
  and set the pointer to the next row after the deleted row. 
  If the deleted row was the last row of the matrix, 
  then the pointer is set to END another time
  and we go up one level further. We also set the pointer to END
  when after a summation at least one entry of the sum vector 
  is larger than two: such sum vectors are \emph{invalid}.

\medskip

  One further case to set the pointer to END is when
  the current sum vector has an entry $1$ at its, say, $k$-th
  position and the current pointer points to a row that
  has missing entries at its positions $1$ to $k$:
  such a sum vector
  can never be completed to a closed vector by adding
  rows that come after the current row. (We use an auxiliary
  variable to keep track of the position of the first non-zero
  entry of the row to which the pointer currently points~to.)

\medskip

  For the above matrix we get the following sequence of combinations:

    \addtolength{\arraycolsep}{1mm} 
  $$\begin{array}{l@{\hspace{4mm}}c@{\hspace{3mm}}c@{\hspace{3mm}}c@{\hspace{3mm}}c@{\hspace{3mm}}c@{\hspace{3mm}}c@{\hspace{4mm}}l}
  -:   & ( & 0 & 0 & 0 & 0 & ) & \mbox{Set pointer to $a$}.\\
  a:   & ( & 2 & 1 & 0 & 0 & ) & \mbox{First entry is $2$: set pointer to $c$.}\\
  a+c: & ( & 2 & 2 & 2 & 0 & ) & \mbox{\emph{Candidate!} Set pointer to END.} \\
  a:   & ( & 2 & 1 & 0 & 0 & ) & \mbox{Set pointer to $d$} \\  
  a+d: & ( & 2 & 2 & 0 & 2 & ) & \mbox{\emph{Candidate!} Set pointer to END.} \\
  a:   & ( & 2 & 1 & 0 & 0 & ) & \mbox{Set pointer to END.} \\  
  -:   & ( & 0 & 0 & 0 & 0 & ) & \mbox{Set pointer to $b$.} \\
  b:   & ( & 2 & 1 & 2 & 2 & ) & \mbox{Set pointer to $c$.} \\
  b+c: & ( & 2 & 2 & 4 & 2 & ) & \mbox{Invalid combination! Set pointer to END.} \\
  b:   & ( & 2 & 1 & 2 & 2 & ) & \mbox{Set pointer to $d$.} \\  
  b+d: & ( & 2 & 2 & 2 & 4 & ) & \mbox{Invalid combination! Set pointer to END.} \\
  b:   & ( & 2 & 1 & 2 & 2 & ) & \mbox{Set pointer to END.} \\  
  -:   & ( & 0 & 0 & 0 & 0 & ) & \mbox{Set pointer to $c$.} \\
  c:   & ( & 0 & 1 & 2 & 0 & ) & \mbox{Set pointer to $d$.} \\
  c+d: & ( & 0 & 2 & 2 & 2 & ) & \mbox{\emph{Candidate!} Set pointer to END.} \\
  c:   & ( & 0 & 1 & 2 & 0 & ) & \mbox{Set pointer to END.} \\  
  -:   & ( & 0 & 0 & 0 & 0 & ) & \mbox{Set pointer to $d$.} \\
  d:   & ( & 0 & 1 & 0 & 2 & ) & \mbox{Set pointer to END.} \\  
  -:   & ( & 0 & 0 & 0 & 0 & ) & \mbox{Set pointer to END.} \\
  \end{array}$$
  \addtolength{\arraycolsep}{-.85mm} 

\medskip

  Thus, for the above example there are three valid combinations,
  $a+c$, $a+d$, and $c+d$, which are further examined 
  in Step~3 and Step~4.

\bigskip

\item[\emph{Step}\hspace{.25mm}] \emph{\hspace{-5.5mm} 3 (Combinatorial Tests):
  Remove complexes from the list of candidates that cannot be manifolds.}

\medskip

  Pure simplicial complexes with the
  pseudomani\-fold property can be regarded as the most general form of pseudomanifolds, 
  as they comprise proper \emph{combinatorial manifolds} (on which we
  will concentrate in the following) as well as \emph{combinatorial pseudomanifolds}, 
  in particular, \emph{Eulerian mani\-folds} 
(see \cite[Ch.~3]{Lutz1999}). 
  For every simplicial complex that we found in Step~2, 
  we perform simple tests to exclude complexes that cannot be
  connected manifolds. 

\medskip

  We first test whether the candidate complex is connected: For \mbox{example},
  the $1$-dimensional complex consisting of the edges $13$, $15$, $24$, $26$, 
  $35$, and $46$ is invariant under the cyclic shift $(1,2,3,4,5,6)$
  and thus is a vertex-transitive pure simplicial complex with the
  pseudomanifold property. However, it is not connected: it is the union of
  two disjoint (empty) triangles.
  In order that a connected simplicial complex is a combinatorial manifold 
  the link of any proper face has to be a combinatorial sphere. 
  Two necessary conditions for this are that the links are connected
  and have the Euler characteristic of a sphere.
  In our implementation, we test these conditions for the link of one vertex $v_0$ (transitivity),
  for the link of every edge containing $v_0$ if $d\geq 3$, and for the link of
  every triangle containing $v_0$ if $d\geq 4$. (If the number of
  vertices is restricted to $n\leq 15$ it is, in most cases, not necessary,
  but expensive to test the links of
  higher-dimensional faces.)

\medskip

  These tests have to be altered only slightly 
  if we want to enumerate all vertex-transitive Eulerian manifolds 
  or all vertex-transitive combinatorial pseudomanifolds 
  for a given vertex-transitive group action;
  see \cite[Ch.~3]{Lutz1999}.

\bigskip

\item[\emph{Step}\hspace{.25mm}] \emph{\hspace{-4mm}4 (Combinatorial Equivalence):
  Remove complexes 
  from the list of candidates that, up to combinatorial equivalence,
  have appeared before.} 

\medskip

  We next classify, up to \emph{combinatorial equivalence} 
  (i.e., up to relabeling the vertices), the candidates which survived Step~3. 
  Two basic combinatorial invariants are particularly helpful for this: 
  the $f$-vector and the \emph{Altshuler-Steinberg
  determinant} \cite{AltshulerSteinberg1973} of a candidate, i.e., the determinant
  ${\rm det}(AA^T)$ of the vertex-facet incidence matrix~$A$ of the
  candidate complex.\index{Altshuler-Steinberg determinant}
  Clearly, ${\rm det}(AA^T)$ is invariant under relabeling vertices or facets.

\medskip

  If the $f$-vectors and the Altshuler-Steinberg determinants of two
  candidate complexes coincide, then one possibility for a
  combinatorial equivalence between these two ver\-tex-tran\-si\-tive complexes 
  is that they are mapped onto each other by an outer automorphism of the acting group. 
  Since many group actions have the cyclic group ${\mathbb Z}_n$, generated by the cycle
  $(1,2,3,\ldots,n)$, as a transitively acting subgroup, we restrict our attention to
  \emph{multiplications}\, $k\mapsto (m\cdot k)\,{\rm mod}\,n$\, with\, 
  $m\in \{ 1,2,3,\ldots,(n-1)\}$\, and\, ${\rm gcd}(m,n)=1$. 

\medskip

  In the above example, the generating simplices
  of the orbits $a$, $b$, $c$, and~$d$ are $1234_{\,7}$, $1235_{\,14}$,
  $1245_{\,7}$, and $1246_{\,7}$, respectively (\emph{the lower index
  indicates the size of the corresponding orbit}). The union of orbits $a+c$\, is 
  mapped to $a+d$\, by multiplication with $2$, and to $c+d$\, by 
  multiplication with $3$. Thus, there is, up to combinatorial equivalence,
  a unique combinatorial $3$-manifold with $7$ vertices and vertex-transitive
  $D_7$-action. This manifold is (of course) the boundary
  complex $\partial C_4(7)$ of the cyclic $4$-polytope $C_4(7)$ with $7$ vertices.

\medskip

  If the $f$-vectors and Altshuler-Steinberg determinants of two
  candidate complexes are equal, but the complexes are not multiplication isomorphic, 
  then we take one simplex of the first complex and test for all possible ways 
  it can be mapped to the generating simplices of the orbits of the second complex
  whether this map can be extended to a simplicial isomorphism of the
  two complexes. (By strong connectivity, a combinatorial isomorphism 
  between two simplicial manifolds is already determined by its action
  on one simplex.)

\medskip

  Alternatively, one can use McKay's (fast!) graph isomorphism testing program
  \texttt{nauty} \cite{nauty} to determine whether the
  vertex-facet incidence graphs of the two complexes are isomorphic or not.
  
\bigskip

\item[\emph{Step}\hspace{.25mm}]
  \emph{\hspace{-5.5mm}s 2\, to\, 4 (Integration of the Combinatorial Steps):
  Whenever we find a \linebreak candidate complex in Step~2,
  we immediately perform Steps~3 to 4 for this complex before we continue with
  the backtracking of Step~2.}

\medskip

  An integration of Steps 2 to 4 has the advantage
  that we do not need to store complexes that would be ruled out by
  Steps~3 to 4
  all along the backtracking. In fact, the basic combinatorial tests of Step~3
  are sufficient to reject most of the candidates that are not
  manifolds (at least for vertex-transitive triangulations with $n\leq 15$ vertices).
  If $d\leq 3$, then all the resulting complexes after Step~3 are
  indeed manifolds: The vertex-links in a combinatorial $2$-mani\-fold are circles,
  whereas the vertex-links in a triangulated $3$-manifold are
  combinatorial $2$-spheres. These are recognized by the
  combinatorial tests of Step~3.  
 
\medskip

  Our GAP-program MANIFOLD\_VT \cite{Lutz_MANIFOLD_VT} is an implementation
  of the Steps 1 to 4 above. All candidate complexes that remain 
  after the backtracking of the integrated Steps~2 to 4 
  (together with the ver\-tex-link of one vertex $v_0$ for each complex if $d\geq 4$)
  are printed to a file. 

\bigskip

\item[\emph{Step}\hspace{.25mm}] \emph{\hspace{-5.5mm} 5 (Homology Computation):
   Remove complexes 
  from the list of candidates for which their homology groups do not
  obey Poincar\'e duality (with respect to\, ${\mathbb Z}_2$-coefficients)
  or for which the homology of the vertex-link differs from the
  homology of a $(d-1)$-sphere.}

\medskip

  We made use of the $C$-program \texttt{homology} by\index{homology}
  Hecken\-bach \cite{Heckenbach1997} to compute the homology groups
  for the candidate complexes and their vertex-links.
  Alternatively, one can use the (proposed) GAP share package
  \emph{Simplicial Homology} \cite{DumasHeckenbachSaundersWelker2003b}
  by Dumas, Heckenbach, Saunders, and Welker 
  (see \cite{DumasHeckenbachSaundersWelker2003a} for a description)
  or the \texttt{TOPAZ} module of the \texttt{polymake} system \cite{polymake}
  of Gawrilow and Joswig to compute the homology groups.
  
\medskip

  In the case of vertex-transitive triangulations with $n\leq 15$ 
  vertices, the candidates that remained after Step~5 
  are all combinatorial manifolds, as we verified in Step~6.

\bigskip

\item[\emph{Step}\hspace{.25mm}] \emph{\hspace{-4.5mm}6 (Recognition of the Vertex-Links):
  Use bistellar flips
  as a heuristics to recognize the vertex-links as combinatorial spheres.} 

\medskip

  We used the program BISTELLAR \cite{Lutz_BISTELLAR}
  (cf.\ also \cite{BjoernerLutz2000}) 
  to check whether the link of a candidate is bistellarly equivalent and
  therefore PL homeo\-morphic to the boundary of a simplex. 
  With this heuristic, it was possible to show that 
  all remaining vertex-transitive candidates with $n\leq 15$ vertices
  are indeed combinatorial manifolds.
  (A fast implementation of the bistellar flip heuristic due to
  Nikolaus Witte is accessible via the \texttt{TOPAZ} module 
  of the \texttt{polymake} system \cite{polymake}.)
  

\bigskip

\item[\emph{Step}\hspace{.25mm}] \emph{\hspace{-4.5mm}7 (Topological Type):
  Use bistellar flips
  or topological classification theorems to determine the topological
  types of the manifolds.} 

\medskip

  For all but one vertex-transitive combinatorial manifold with $n\leq 15$ 
  vertices that we found in the previous steps it has been possible to
  determine their homeomorphism types. 
  This was done in most cases with the program BISTELLAR\_EQUIVALENT \cite{Lutz_BISTELLAR_EQUIVALENT},
  which we used to establish a bistellar equivalence between the test manifold  
  and some reference mani\-fold. As reference manifolds we used small or minimal triangulations 
  of manifolds that have the same homology as the test object, but for
  which their topological types are known. In a small number of cases,
  topological classification theorems were used to determine the 
  topological types. For details see Section~\ref{sec:top_types}.
 
\medskip

  Figure~\ref{fig:torus-flips} displays an application of bistellar flips: 
  We explicitly give a bistellar equivalence between a $9$-vertex\index{torus!M\"obius}
  triangulation of the $2$-torus and the unique vertex-mini\-mal
  $7$-vertex triangulation of M\"obius~\cite{Moebius1886}.

\bigskip

  \begin{figure}
    \begin{center}
      \psfrag{9-vertex torus}{9-vertex torus}
      \psfrag{Mobius' 7-vertex torus}{M\"obius' 7-vertex torus}
      \psfrag{.}{}
      \includegraphics[width=.85\linewidth]{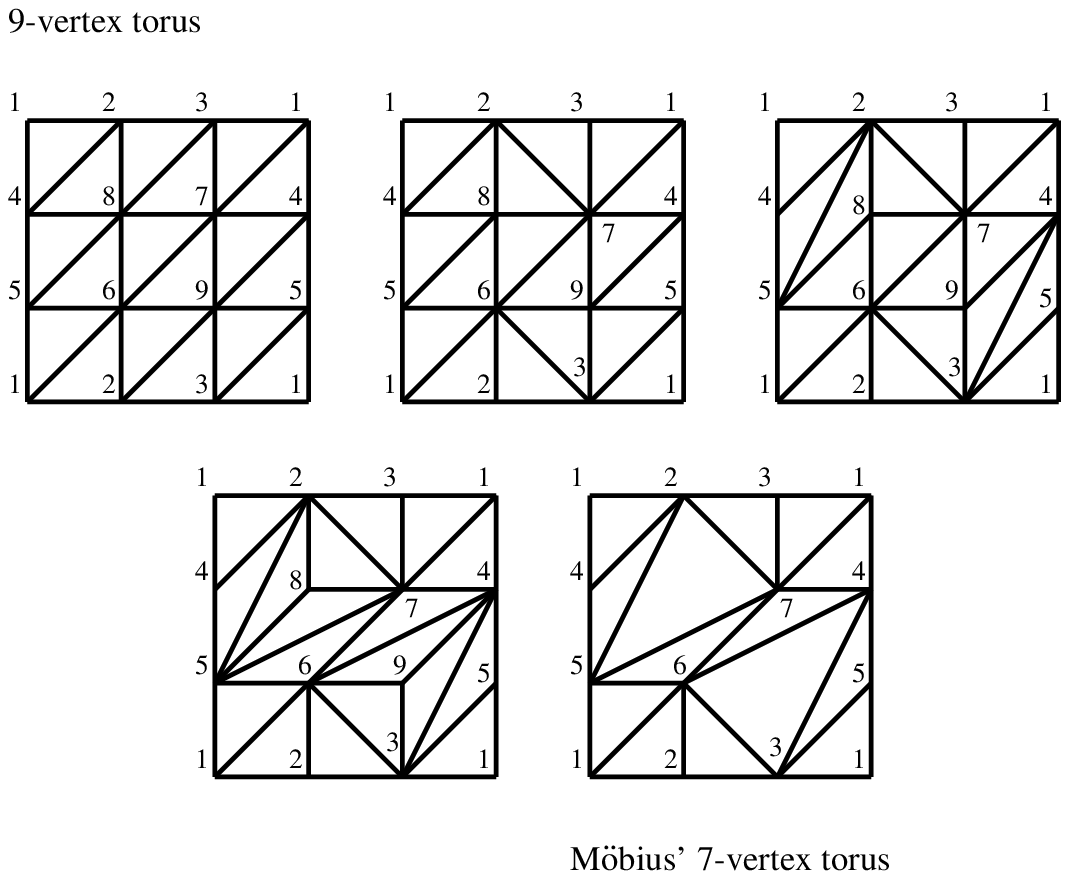}
    \end{center}
    \caption{Bistellar flips on the $9$-vertex torus to reduce the number of vertices.}
    \label{fig:torus-flips}
  \end{figure}

\end{list}

\noindent
With respect to computation time we remark that the most expensive
step in the above procedure is the backtracking of Step~2.
More precisely, the computation time crucially depends on
the size of the associated matrices. As a consequence,
vertex-transitive triangulations for group actions
of small size are hard to enumerate. If, on the other side, 
the acting group has large group order, then there is a good chance
to complete the enumeration -- even for larger $n$ and $d$.

The division line between instances that can be computed
and those that cannot is sharp: It takes minutes or at most hours
to enumerate vertex-transitive triangulations with a dihedral
group action on $14$ and $15$ vertices, but it is hopeless
to complete the enumeration (at least with the present techniques)
for cyclic actions on $14$ and $15$ vertices in dimensions $5\leq d\leq7$.


\section{Enumeration Results}

We used the program MANIFOLD\_VT \cite{Lutz_MANIFOLD_VT} to enumerate 
(candidates for) ver\-tex-transitive triangulated $d$-manifolds
with $n\leq 15$ vertices and $2\leq d\leq n-2$ for almost
all the actions of transitive permutation groups of degree $n\leq 15$.
The cases where an enumeration was not possible are
for $5\leq d\leq 7$
the actions 
$\mbox{}^5\hspace{.3pt}14^{\,1}$, 
$\mbox{}^6\hspace{.3pt}14^{\,1}$, 
$\mbox{}^7\hspace{.3pt}14^{\,1}$
of the cyclic group ${\mathbb Z}_{14}$ 
and the actions
$\mbox{}^5\hspace{.3pt}14^{\,2}$,
$\mbox{}^6\hspace{.3pt}14^{\,2}$, 
$\mbox{}^7\hspace{.3pt}14^{\,2}$
of the dihedral group $D_7$ on $14$ vertices
as well as for $4\leq d\leq 8$
the actions
$\mbox{}^4\hspace{.3pt}15^{\,1}$,
$\mbox{}^5\hspace{.3pt}15^{\,1}$,
$\mbox{}^6\hspace{.3pt}15^{\,1}$,
$\mbox{}^7\hspace{.3pt}15^{\,1}$, 
$\mbox{}^8\hspace{.3pt}15^{\,1}$
of the cyclic group ${\mathbb Z}_{15}$ on $15$ vertices.
For all but these $11$ actions we succeeded
to complete the enumeration (Steps 2 to 4) and also Steps~5 and 6 
of the previous section.
All candidates that remained after Step~6 
turned out to be combinatorial manifolds.

\begin{thm} 
There are at least\, $525$ combinatorial manifolds of dimension\, 
$2\leq d\leq 13$\, with\, $n\leq 15$\, vertices 
that have a vertex-transitive automorphism group. 
\end{thm}
According to the Brehm-K\"uhnel bound~\cite{BrehmKuehnel1987},
every combinatorial manifold with $n < 3\left\lceil\frac{d}{2}\right\rceil +3$
vertices is a sphere and every combinatorial manifold with
$n= 3\left\lceil\frac{d}{2}\right\rceil +3$ vertices
is either a sphere or a manifold `like a projective plane'.
The latter case is only possible for $d=2,4,8$, and $16$ 
(cf.\ \cite{BrehmKuehnel1987} and also 
\cite{Lutz2005bpre}).
Among the $525$ mani\-folds that we found there are $235$ spheres
and $290$ examples of other topological types. The 
explicit numbers of spheres and non-spheres are listed 
for given $d$ and $n$ (together with the Brehm-K\"uhnel bound for $2\leq d\leq 8$)
in Table~\ref{tbl:main-table}.

\begin{cor}
There are precisely\, $220$ combinatorial manifolds of dimension\,
$2\leq d\leq 11$\, with\, $n\leq 13$ vertices 
that have a vertex-transitive automorphism group:
$110$ spheres and $110$ manifolds that are not spheres.
The\, $34$ different homeomorphism types of these manifolds are for

\medskip

\begin{tabular}{l@{\hspace{3mm}}l}
$d=2$: & $S^2$, ${\bf T}^2$, the orientable surfaces of genus\, $2$, $3$, $4$, $5$, $6$,\\[1.5mm]
       & ${\mathbb R}{\bf P}^{\,2}$, the non-orientable surfaces of
          genus\, $2$, $4$, $5$, $7$, $8$, $15$,\\[1.5mm]
$d=3$: & $S^3$, $S^2\hbox{$\times\hspace{-1.62ex}\_\hspace{-.4ex}\_\hspace{.7ex}$}S^1$,
         $S^2\!\times\!S^1$, $(S^2\!\times\!S^1)^{\# 2}$, ${\mathbb R}{\bf P}^{\,3}$,\\[1.5mm]
$d=4$: & $S^4$, ${\mathbb C}{\bf P}^{\,2}$, $S^3\!\times\!S^1$, 
         $S^3\hbox{$\times\hspace{-1.62ex}\_\hspace{-.4ex}\_\hspace{.7ex}$}S^1$, 
         $S^2\!\times\!S^2$, $(S^2\!\times\!S^2)^{\# 2}$,\\[1.5mm]
$d=5$: & $S^5$, $S^4\hbox{$\times\hspace{-1.62ex}\_\hspace{-.4ex}\_\hspace{.7ex}$}S^1$, $SU(3)/SO(3)$,
\end{tabular}

\medskip

\noindent
as well as $S^6$, $S^7$, $S^8$, $S^9$, $S^{10}$, and\, $S^{11}$.     
\end{cor}
The topological types of the examples were determined according to
Step~7 of the previous section. For details see Section~\ref{sec:top_types}.

\begin{cor}
There are precisely\, $77$ combinatorial\, $2$-manifolds with $n\leq 15$ vertices 
that have a vertex-transitive automorphism group.
Of these examples $42$ are orientable and\, $35$ are non-orientable.
In particular, there are $18$ different topological types:
$S^2$, ${\bf T}^2$, the orientable surfaces of genus $2$, $3$, $4$, $5$, $6$, $8$, 
${\mathbb R}{\bf P}^{\,2}$, and the non-orientable surfaces of genus
$2$, $4$, $5$, $7$, $8$, $12$, $15$, $16$, and~$17$. 
\end{cor}

\begin{cor}
There are exactly\, $166$ combinatorial\, $3$-manifolds on\, $n\leq 15$ vertices 
with a vertex-transitive automorphism group; $52$ of these are spheres,
whereas $114$ are not spheres. The manifolds are of one of $8$ different topo\-lo\-gical
types: $S^3$, $S^2\hbox{$\times\hspace{-1.62ex}\_\hspace{-.4ex}\_\hspace{.7ex}$}S^1$,
$S^2\!\times\!S^1$, $(S^2\!\times\!S^1)^{\# 2}$, ${\mathbb R}{\bf P}^{\,3}$,
$L(3,1)$, $S^3/Q$, and\, ${\bf T}^3$.
\end{cor}

We label every vertex-transitive combinatorial manifold that we found 
(up to combinatorial equivalence) by our enumeration
with a \emph{unique symbol}: The \mbox{$k$-th} example of a combinatorial manifold 
of dimension $d$ with $n$ vertices that is listed for the $i$-th 
transitive permutation group $n^{\,i}$ of degree $n$ is denoted 
by\, $\mbox{}^d\hspace{.3pt}n^{\,i}_{\,k}$.
We list the respective manifolds in the Tables~\ref{tbl:2man} to \ref{tbl:11man},
together with additional information on their topological types 
and where they appeared previously in the literature -- as far as we know.

\begin{cor}
There is no vertex-transitive triangulation of a combinatorial
$5$-manifold, different from the $5$-sphere, with $12$ vertices.
Also there is no vertex-transitive triangulation of a combinatorial
$6$-manifold, different from the $6$-sphere, with $13$ vertices.
\end{cor}

\begin{landscape}

\begin{table}
\caption{\protect\parbox[t]{6.05in}{Vertex-transitive combinatorial manifolds with $n\leq 15$ vertices; numbers of spheres and non-spheres (in bold).}}
\label{tbl:main-table}
\small\centering
\defaultaddspace=0.3em
\begin{tabular*}{\linewidth}{@{\extracolsep{\fill}}@{}rrrrrrrrrrrrr@{}}   
\addlinespace
\addlinespace
\addlinespace
\toprule
\addlinespace
\addlinespace
$n\backslash d$ &        2 &        3 &            4 &            5 &           6 &            7 &         8 &   9 &  10 &  11 &  12 &  13 \\ [.5mm]\midrule
 \addlinespace
 \addlinespace
     15 &        0 &        5 &     $\geq 2$ &    $\geq 23$ &    $\geq 0$ &    $\geq 14$ &  $\geq 0$ &   5 &   0 &   3 &   0 &   1 \\
        & $\bf 17$ & $\bf 29$ & $\bf\geq 31$ & $\bf\geq 12$ & $\bf\geq 10$&  $\bf\geq 0$ &   $\bf 1$ &     &     &     &     &     \\ [1.5mm]\cline{7-8}
 \addlinespace
     14 &        0 &       20 &            0 &    $\geq 16$ &    $\geq 3$ &     $\geq 6$ &        23 &   1 &   1 &   1 &   1 &     \\
        & $\bf 13$ & $\bf 35$ &     $\bf 23$ &  $\bf\geq 9$ & $\bf\geq 0$ &              &           &     &     &     &     &     \\ 
 \addlinespace
     13 &        0 &        6 &            0 &           17 &           0 &            6 &         0 &   1 &   0 &   1 &     &     \\
        &  $\bf 4$ &  $\bf 9$ &      $\bf 5$ &      $\bf 2$ &     $\bf 0$ &              &           &     &     &     &     &     \\ [1.5mm]\cline{6-6}
 \addlinespace
     12 &        1 &        6 &            1 &           27 &           0 &            4 &         1 &   1 &   1 &     &     &     \\
        & $\bf 30$ & $\bf 33$ &      $\bf 7$ &      $\bf 0$ &             &              &           &     &     &     &     &     \\ [1.5mm]\cline{5-5}
 \addlinespace
     11 &        0 &        3 &            0 &            3 &           0 &            1 &         0 &   1 &     &     &     &     \\
        &  $\bf 1$ &  $\bf 3$ &      $\bf 1$ &              &             &              &           &     &     &     &     &     \\
 \addlinespace
     10 &        0 &        6 &            4 &            1 &           1 &            1 &         1 &     &     &     &     &     \\
        &  $\bf 3$ &  $\bf 4$ &      $\bf 0$ &              &             &              &           &     &     &     &     &     \\
 \addlinespace
      9 &        0 &        1 &            0 &            2 &           0 &            1 &           &     &     &     &     &     \\
        &  $\bf 3$ &  $\bf 1$ &      $\bf 1$ &              &             &              &           &     &     &     &     &     \\ [1.5mm]\cline{3-4}
 \addlinespace
      8 &        0 &        2 &            0 &            1 &           1 &              &           &     &     &     &     &     \\
        &  $\bf 1$ &          &              &              &             &              &           &     &     &     &     &     \\
 \addlinespace
      7 &        0 &        1 &            0 &            1 &             &              &           &     &     &     &     &     \\
        &  $\bf 1$ &          &              &              &             &              &           &     &     &     &     &     \\
 \addlinespace
      6 &        1 &        1 &            1 &              &             &              &           &     &     &     &     &     \\
        &  $\bf 1$ &          &              &              &             &              &           &     &     &     &     &     \\ [1.5mm]\cline{2-2}
 \addlinespace
      5 &        0 &        1 &              &              &             &              &           &     &     &     &     &     \\
\\
 \addlinespace
      4 &        1 &          &              &              &             &              &           &     &     &     &     &     \\
 \addlinespace
 \addlinespace
\bottomrule
\end{tabular*}
\end{table}

\end{landscape}

Nonetheless, there are asymmetric triangulations of $S^3\!\times\!S^2$
with $12$ vertices and of $S^3\!\times\!S^3$ with $13$ vertices
from which it follows that the Brehm-K\"uhnel bound is sharp for $d=5$\, and\, $d=6$,
respectively (see 
\cite{Lutz2005bpre}).
For $d=2,3,4,8$, 
the $6$-vertex triangulation\, $\mbox{}^2\hspace{.3pt}6^{\,12}_{\,1}$\, of the real projective plane,
Walkup's \cite{Walkup1970} $9$-vertex triangulation\,
$\mbox{}^3\hspace{.3pt}9^{\,3}_{\,2}$\, of\, $S^2\hbox{$\times\hspace{-1.62ex}\_\hspace{-.4ex}\_\hspace{.7ex}$}S^1$,
K\"uhnel's \cite{KuehnelBanchoff1983} $9$-vertex triangulation\,
$\mbox{}^4\hspace{.3pt}9^{\,13}_{\,1}$\, of\, ${\mathbb C}{\bf P}^{\,2}$, and
Brehm and K\"uhnel's $A_5$-invariant triangulation\,
$\mbox{}^8\hspace{.3pt}15^{\,5}_{\,1}$\, of a manifold ${\sim}{\mathbb H}{\bf P}^{\,2}$
like the quaternionic projective plane with $15$ vertices are
vertex-transitive examples of combinatorial manifolds that are
vertex-minimal by the Brehm-K\"uhnel bound in the respective dimensions.

\begin{thm}
There are two vertex-minimal, vertex-transitive triangulations\linebreak
$\mbox{}^4\hspace{.3pt}12^{\,2}_{\,1}$ and\,
$\mbox{}^4\hspace{.3pt}12^{\,2}_{\,2}$\, of\,
$(S^2\!\times\!S^2)\# (S^2\!\times\!S^2)$ with $12$ vertices.
\end{thm}
\noindent
\textbf{Proof.} 
Vertex-minimality for these two triangulations follows from
K\"uhnel's bound \cite[4.1]{Kuehnel1990-few}   
which states that $\binom{n-4}{3}\geq 10\,(\chi (M)-2)$ for every
combinatorial $4$-mani\-fold~$M$ with $n$ vertices.
\hfill $\Box$

\section{Topological Types}
\label{sec:top_types}

Various of the $525$ vertex-transitive combinatorial manifolds
that we found with $n\leq 15$ vertices appeared previously
in the literature. For many of these examples, their
topological types were determined in the respective papers;
see the references cited in the Tables~\ref{tbl:3man} to \ref{tbl:11man}.
For all but one of the examples their types can also be recognized
as described in the following.

\medskip

Since the topological type of a $2$-dimensional manifold can be
determined from its Euler characteristic and its
orientation character, the actual work for recognizing topological types
starts with dimension $3$: All vertex-transitive $3$-mani\-folds 
with $n\leq 15$ vertices from our enumeration turned out to be either 
triangulations of Seifert manifolds or triangulations of the connected sum $(S^2\!\times\!S^1)^{\# 2}$.
Reference triangulations for Seifert manifolds are available
via the program SEIFERT~\cite{Lutz_SEIFERT} (for a description 
of the program see \cite{BrehmLutz2002pre} and 
\cite{Lutz2003bpre}).
The connected sum $(S^2\!\times\!S^1)^{\# 2}$ can be composed
combinatorially by taking two disjoint copies of a triangulation of
$S^2\!\times\!S^1$, then removing a simplex each,
and finally glueing both parts together. By using these reference triangulations,
it was possible to recognize the topological types of all the vertex-transitive
$3$-manifolds with $n\leq 15$ vertices with the bistellar flip program
BISTELLAR\_EQUIVALENT~\cite{Lutz_BISTELLAR_EQUIVALENT}.

\medskip

The $4$-dimensional combinatorial manifolds that we found for $n\leq 15$ with a
vertex-transitive group action are of the topological types $S^4$, ${\mathbb C}{\bf P}^{\,2}$, $S^3\!\times\!S^1$, 
$S^3\hbox{$\times\hspace{-1.62ex}\_\hspace{-.4ex}\_\hspace{.7ex}$}S^1$, 
$S^2\!\times\!S^2$, $(S^2\!\times\!S^2)^{\# 2}$, and
$(S^3\hbox{$\times\hspace{-1.62ex}\_\hspace{-.4ex}\_\hspace{.7ex}$}S^1)\,\#\,({\mathbb C}{\bf P}^{\,2})^{\# 5}$.
Reference triangulations for the (twisted) sphere products $S^3\!\times\!S^1$, $S^2\!\times\!S^2$,
and $S^3\hbox{$\times\hspace{-1.62ex}\_\hspace{-.4ex}\_\hspace{.7ex}$}S^1$
can be obtained as (twisted) product triangulations; see 
\cite{Lutz2003bpre}.
As a reference triangulation for $S^4$ we can take the boundary of the
$5$-simplex, and the connected sum \linebreak $(S^2\!\times\!S^2)^{\# 2}$ can be
composed as described above. Some of the vertex-transitive triangulations of
$4$-manifolds of these topological types were known before.
However, it was also possible to recognize all these examples 
with the program BISTELLAR\_EQUIVALENT. The $9$-vertex triangulation 
of K\"uhnel of ${\mathbb C}{\bf P}^{\,2}$ is discussed in \cite{KuehnelBanchoff1983}.
Thus it remained to determine the type of the example $\mbox{}^4\hspace{.3pt}15^{\,4}_{\,1}$.

\begin{thm}
There is a vertex-transitive triangulation $\mbox{}^4\hspace{.3pt}15^{\,4}_{\,1}$
of the manifold 
$(S^3\hbox{$\times\hspace{-1.62ex}\_\hspace{-.4ex}\_\hspace{.7ex}$}S^1)\#\,({\mathbb C}{\bf P}^{\,2})^{\# 5}$
with $15$ vertices.
\end{thm}
\noindent
\textbf{Proof.} 
The triangulation $\mbox{}^4\hspace{.3pt}15^{\,4}_{\,1}$ has 
infinite cyclic fundamental group, which we computed with the 
group algebra package GAP~\cite{GAP4}. 
By the classification of Wang \cite{Wang1995} of closed non-orientable $4$-manifolds 
with infinite cyclic fundamental group,
$(S^3\hbox{$\times\hspace{-1.62ex}\_\hspace{-.4ex}\_\hspace{.7ex}$}S^1)\#\,({\mathbb C}{\bf P}^{\,2})^{\# 5}$
is the only such $4$-manifold with homology $H_*=({\mathbb Z},{\mathbb Z},{\mathbb Z}^5,{\mathbb Z}_2,0)$.
\hfill $\Box$

\begin{conj}
The combinatorial manifold $\mbox{}^4\hspace{.3pt}15^{\,4}_{\,1}$
is the unique vertex-mini\-mal triangulation of 
$(S^3\hbox{$\times\hspace{-1.62ex}\_\hspace{-.4ex}\_\hspace{.7ex}$}S^1)\#\,({\mathbb C}{\bf P}^{\,2})^{\# 5}$
with $15$ vertices.
\end{conj}

The vertex-transitive combinatorial $5$-manifolds that we found with
$n\leq 15$ vertices are of the topological types\,
$S^5$, $S^4\hbox{$\times\hspace{-1.62ex}\_\hspace{-.4ex}\_\hspace{.7ex}$}S^1$, $SU(3)/SO(3)$,
and $S^3\!\times\!S^2$. All the corresponding triangulations of
$S^5$ and $S^4\hbox{$\times\hspace{-1.62ex}\_\hspace{-.4ex}\_\hspace{.7ex}$}S^1$
could be recognized with bistellar flips.

\begin{thm}
There is a $3$-neighborly triangulation $\mbox{}^5\hspace{.3pt}13^{\,3}_{\,2}$\, 
of the simply connected homogeneous $5$-dimensional manifold\, $SU(3)/SO(3)$
with homology $H_*=({\mathbb Z},0,{\mathbb Z}_2,0,0,{\mathbb Z})$.
This triangulation 
with $f=(13,\underline{78},\underline{286},533,468,156)$ 
has a vertex-transitive action of the affine group $13\!:\!3$.
\end{thm}
\noindent
\textbf{Proof.} 
Since the triangulation $\mbox{}^5\hspace{.3pt}13^{\,3}_{\,2}$\, is
$3$-neighborly, the corresponding mani\-fold is simply connected.
According to the classification of all simply connected $5$-manifolds
by Barden~\cite{Barden1965}, there is only one simply connected
$5$-manifold with homology $H_*=({\mathbb Z},0,{\mathbb Z}_2,0,0,{\mathbb Z})$, which he denoted by $X_{-1}$.
In fact, it is the well-known simply connected homogeneous $5$-manifold $SU(3)/SO(3)$;
cf.\ the classification of compact homogeneous manifolds of low\index{manifold!aaSU(3)/SO(3)@$SU(3)/SO(3)$}
dimension by Gorbatsevich~\cite{Gorbatsevich1980} as well as the
exposition on homogeneous manifolds in \cite{Klaus1988}.
\mbox{}\hfill $\Box$              

\begin{conj}
The $3$-neighborly triangulation $\mbox{}^5\hspace{.3pt}13^{\,3}_{\,2}$ 
is the unique vertex-mini\-mal triangulation of $SU(3)/SO(3)$
with $13$ vertices.
\end{conj}
There are two further vertex-transitive triangulations 
$\mbox{}^5\hspace{.3pt}14^{\,4}_{\,2}$ and
$\mbox{}^5\hspace{.3pt}14^{\,4}_{\,6}$ of \linebreak 
$SU(3)/SO(3)$ that are bistellarly equivalent to the triangulation
$\mbox{}^5\hspace{.3pt}13^{\,3}_{\,2}$. However, we did not find another
triangulation of $SU(3)/SO(3)$ with $13$ vertices.

\medskip

Altogether, our enumeration yielded four vertex-transitive triangulations
$\mbox{}^5\hspace{.3pt}14^{\,3}_{\,8}$, 
$\mbox{}^5\hspace{.3pt}14^{\,3}_{\,13}$,
$\mbox{}^5\hspace{.3pt}14^{\,3}_{\,14}$, and 
$\mbox{}^5\hspace{.3pt}14^{\,3}_{\,15}$
of $S^3\!\times\!S^2$ with $n\leq 15$ vertices:\index{manifold!aaS3xS2@$S^3\times S^2$}
For all four examples we used the program BISTELLAR~\cite{Lutz_BISTELLAR}
to obtain vertex-minimal $3$-neighborly triangulations 
with $12$ vertices and $f=(12,\underline{66},\underline{220},390,336,112)$; see 
\cite{Lutz2005bpre}.
In particular, it follows that the four examples are simply-connected.
By the classification of Barden~\cite{Barden1965}, there are precisely
two simply connected $5$-manifolds with the homology of $S^3\!\times\!S^2$,
namely $M_{\infty}=S^3\!\times\!S^2$ with trivial
and $X_{\infty}$ with non-vanishing second Stiefel-Whitney class.
The first example, $\mbox{}^5\hspace{.3pt}14^{\,3}_{\,8}$,
is centrally symmetric (cf.\ 
\cite{Lutz2004apre})
and is therefore embedded in the $6$-di\-men\-sional boundary complex 
$\partial C_7^{\Delta}$ of the $7$-dimensional crosspolytope
$C_7^{\Delta}$ with $14$ vertices.
Since $\mbox{}^5\hspace{.3pt}14^{\,3}_{\,8}$ is a codimension
$1$-submanifold in the sphere $\partial C_7^{\Delta}$,
it divides $\partial C_7^{\Delta}$ into two connected components
by Alexander duality, 
both parts having $\mbox{}^5\hspace{.3pt}14^{\,3}_{\,8}$
as their common boundary. By a theorem of Pontrjagin,   
the Stiefel-Whitney numbers of a $d$-manifold that is the boundary of a
smooth compact $(d+1)$-manifold are all zero (cf.~\cite[4.9]{MilnorStasheff1974}).
Hence, $\mbox{}^5\hspace{.3pt}14^{\,3}_{\,8}$ is a triangulation
of $S^3\!\times\!S^2$. For the other three triangulations
$\mbox{}^5\hspace{.3pt}14^{\,3}_{\,13}$,
$\mbox{}^5\hspace{.3pt}14^{\,3}_{\,14}$, and 
$\mbox{}^5\hspace{.3pt}14^{\,3}_{\,15}$ we computed their 
Stiefel-Whitney classes with the \texttt{TOPAZ} module of the
\texttt{polymake} system \cite{polymake}.
In all three cases the Stiefel-Whitney classes vanish
and the examples are therefore triangulations of $S^3\!\times\!S^2$.

\medskip

The vertex-transitive combinatorial $6$-manifolds with $n\leq 15$
vertices that we obtained by our enumeration are of three topological
types, $S^6$, $S^5\!\times\!S^1$, and $S^3\!\times\!S^3$.
The triangulations of $S^6$ were recognized with bistellar flips and 
the one triangulation $\mbox{}^6\hspace{.3pt}15^{\,2}_{\,1}$
of\, $S^5\!\times\!S^1$ is a member of a series of 
sphere products of K\"uhnel~\cite[$M^6$]{Kuehnel1986a-series}
(see also \cite[$M^6_5$]{KuehnelLassmann1996-bundle}).

For the remaining vertex-transitive combinatorial $6$-manifolds
we made use of the classification of $2$-connected topological 
$6$-manifolds by ${\rm \check{Z}}$ubr~\cite{Zubr1988}
(cf.\ also \cite{Kreck2001} and 
\cite{KuehnelLutz2003pre}):
The topological type of a $2$-connected $6$-manifold is determined by its
Euler characteristic. We used two approaches to show that the 
examples are $2$-connected. In both approaches, we first
computed the homology groups of the examples, which are the homology
groups of $S^3\!\times\!S^3$. Then, in the first approach, 
we computed the fundamental group of the examples with the program GAP, 
which turned out to be trivial in each case. 
According to the theorem of Hurewicz (cf.\ \cite[p.~80]{Matousek2003}), 
every simply connected space with trivial second homology is $2$-connected. 
Thus, by the classification of ${\rm \check{Z}}$ubr it follows that
the examples are triangulations of $S^3\!\times\!S^3$.
Another way to show the $2$-connectedness is by using
bistellar flips to reduces the examples to
vertex-minimal triangulations with $13$ vertices.
The resulting $f$-vector that we obtained in each of the cases is
$f=(13,\underline{78},\underline{286},\underline{715},1014,728,208)$.
From the $f$-vector it can be read of that these $13$-vertex 
triangulations are $4$-neighborly and therefore $2$-connected.

\medskip

Apart from the combinatorial manifold $\mbox{}^8\hspace{.3pt}15^{\,5}_{\,1}$,
all vertex-transitive triangulations of manifolds of dimension $7\leq d\leq 13$
that we found with $n\leq 15$ vertices are spheres, as we recognized
with bistellar flips.

The vertex-transitive example $\mbox{}^8\hspace{.3pt}15^{\,5}_{\,1}$ was first
discovered by Brehm and K\"uhnel \cite[$M^8_{15}$]{BrehmKuehnel1992}.
It has homo\-logy $H_*=({\mathbb Z},0,0,0,{\mathbb Z},0,0,0,{\mathbb Z})$.
According to the Brehm-K\"uhnel bound \cite{BrehmKuehnel1987},
every combinatorial $8$-manifold with $15$ vertices is either a sphere
or a mani\-fold like the quaternionic projective plane\index{projective space!`quaternionic'}
(in the sense of\, {\rm\cite{EellsKuiper1962b}}).
There are infinitely many such manifolds that
can be distinguished by their first Pontrjagin class.  
However, it is unclear how to explicitely compute the first Pontrjagin class
combinatorially for simplicial complexes of the size of $\mbox{}^8\hspace{.3pt}15^{\,5}_{\,1}$.
We denote the topological type of the example $\mbox{}^8\hspace{.3pt}15^{\,5}_{\,1}$
by ${\sim}{\mathbb H}{\bf P}^{\,2}$ to indicate that, most likely,
it is homeomorphic to ${\mathbb H}{\bf P}^{\,2}$.

With the exception of $8$-manifolds that possibly have the cyclic
permutation group $15^{\,1}$
as vertex-transitive automorphism group, we were able to enumerate all 
vertex-transitive $8$-manifolds with $n\leq 15$ vertices (cf.\ Table~\ref{tbl:8man}).
Besides $\mbox{}^8\hspace{.3pt}15^{\,5}_{\,1}$, all the respective examples are spheres.
The combinatorial mani\-fold $\mbox{}^8\hspace{.3pt}15^{\,5}_{\,1}$ of Brehm
and K\"uhnel has the group $A_5$ as its vertex-transitive 
automorphism group and is the only example with this group action. 

\begin{cor}  {\rm (Brehm~\cite{Brehm-pers})}
There is exactly one vertex-transitive triangulation
of a manifold like the quaternionic projective plane
with $15$ vertices.
\end{cor}

\noindent
\textbf{Proof.} 
Manifolds like the quaternionic projective plane are $3$-connected 
and have Euler characteristic $\chi =3$. By \cite[4.7]{Kuehnel1995-book}
it follows that a triangulation of a manifold like the quaternionic projective plane
with $15$ vertices is $5$-neigh\-borly. By the Dehn-Sommerville equations
the Euler characteristic of a $5$-neighborly combinatorial $8$-manifold 
completely determines the $f$-vector  (see the discussion in
\cite{BrehmKuehnel1992} and in the proof of Theorem~4.17 of \cite{Kuehnel1995-book}).
In the case of $n=15$ vertices and Euler characteristic $\chi =3$ the resulting $f$-vector is 
$f=(15,\underline{105},\underline{455},\underline{1365},\underline{3003},4515,4230,2205,490).$

It remains to rule out that there are triangulations of manifolds like the quaternionic projective plane
with $15$ vertices and the above $f$-vector that are invariant under
the action of the cyclic group $15^{\,1}$ with generator $(1,2,3,\dots,15)$.
This was done by Brehm~\cite{Brehm-pers}: The cyclic action $15^{\,1}$ has $335$ orbits 
of $9$-tuples and, by complementarity (cf.~\cite[p~75]{Kuehnel1995-book}), 
also $335$ orbits of $6$-tuples, $333$ of size $15$ and $2$ of size $5$ in both cases. 
In order to compose an $8$-manifold with $490$ $8$-simplices,
both small orbits of $9$-tuples of size $5$ have to be used.
These orbits are generated by the simplices
$123678\,11\,12\,13$ and $124679\,11\,12\,14$.
On the other hand, the two orbits of $6$-tuples of size $5$
cannot be taken as $5$-faces since $f_5=4515$. But these two orbits
are generated by the simplices $1267\,11\,12$ and $1368\,11\,13$
and are thus included as subfaces in the above orbits. Contradiction. 
\hfill $\Box$

\section{Tables of Manifolds}
\label{sec:tables_manif}

We list the combinatorial manifolds that we found with a vertex-transitive group action 
on $n\leq 15$ vertices in the Tables~\ref{tbl:2man} to \ref{tbl:11man}.
For every example we give the lexicographically smallest 
simplices of the respective orbits of facets as orbit
representatives. 
All examples can be rebuilt from their orbit
representatives by using GAP commands as follows:

\pagebreak

\begin{verbatim}
gap> G:=TransitiveGroup(7,4);
gap> facets:=[ ];
gap> UniteSet(facets,Orbit(G,[1,2,4],OnSets));
gap> Print(facets,"\n");
[ [ 1, 2, 4 ], [ 1, 2, 6 ], [ 1, 3, 4 ], [ 1, 3, 7 ],
  [ 1, 5, 6 ], [ 1, 5, 7 ], [ 2, 3, 5 ], [ 2, 3, 7 ],
  [ 2, 4, 5 ], [ 2, 6, 7 ], [ 3, 4, 6 ], [ 3, 5, 6 ],
  [ 4, 5, 7 ], [ 4, 6, 7 ] ]
\end{verbatim}
\medskip

\noindent
This example is M\"obius' $7$-vertex torus $\mbox{}^2\hspace{.3pt}7^{\,4}_{\,1}$
with vertex-transitive automorphism group $7^{\,4}$.
The triangle $124$ (we omit brackets and commas in the tables to save space)
is listed as the generating triangle for the orbit, since it
is the lexicographically smallest triangle in the orbit.
Examples with more than one
orbit of facets can be built similarly by uniting their orbits
of facets with additional\,\, \texttt{UniteSet}\,\, commands.
The lower index attached to every generating simplex in the tables indicates
the size of the corresponding orbit. For M\"obius' torus
the orbit $124_{\,14}$ has the above $14$ triangles.

\medskip

For every transitive permutation group of degree $n\leq 15$ 
that occurs as the automorphism group of one of the combinatorial
manifolds in the Tables~\ref{tbl:2man} to \ref{tbl:11man},
we list the generators of the group in Table~\ref{tbl:generators},
with the following exceptions: 
The cyclic group actions, $10^1$, $11^1$, $12^1$, $13^1$, $14^1$, and $15^1$,  
and the dihedral group actions, $7^2$, $8^6$, $9^3$, $10^3$, $11^2$, $12^{12}$, $13^2$, $14^3$, and $15^2$,
with generators $a_n=(1,2,3,\ldots ,n)$ and $b_n=(1,2\lfloor \frac{n}{2}\rfloor )(2,2\lfloor \frac{n}{2}\rfloor\!-\!1)\ldots (\lfloor \frac{n}{2}\rfloor,\lfloor \frac{n}{2}\rfloor\!+\!1)$,
such that ${\mathbb Z}_n=\langle a_n\rangle$ and $D_n=\langle a_n,b_n\rangle$,
are not listed in Table~\ref{tbl:generators}.
Also the automorphism groups of the boundary complexes of
$(d+1)$-simplices, i.e., the symmetric groups $S_{d+2}$ of $d+2$
elements, where $2\leq d\leq 13$, are omitted from Table~\ref{tbl:generators}.

\begin{cor}
The transitive permutation groups of Table~\mbox{\rm\ref{tbl:generators}}
together with the additional groups above are precisely all permutation groups 
that occur as vertex-transitive automorphism groups of combinatorial manifolds 
with\, $n\leq 15$ vertices.
\end{cor}

Kimmerle and Kouzoudi \cite{KimmerleKouzoudi2003} determined
that the boundary complex of the $3$-simplex with $4$ vertices, 
the real projective plane $\mbox{}^2\hspace{.3pt}6^{\,11}_{\,1}$,
triangulated minimally with $6$ vertices, and
M\"obius' vertex-minimal $7$-vertex torus $\mbox{}^2\hspace{.3pt}7^{\,4}_{\,1}$
are the only combinatorial surfaces that admit a \emph{doubly transitive}
automorphism group.

Not all of the transitive group actions from Table~\ref{tbl:generators} have systematic names.
We use the GAP terminology for these groups:
e.g., $t8n15(32)$\, is the transitive permutation group
on $8$ vertices with number $15$ (we add the size of the group in brackets).
There are finite groups that have more than one representation 
as a transitive permutation group on $n$ vertices. For example,
the permutation groups $t12n8(24)$ and $t12n9(24)$ are 
different representations of the symmetric group $S_4$.

\pagebreak

Further symbols and abbreviations that are frequently used 
in the Tables~\ref{tbl:2man} to \ref{tbl:11man} are:

\medskip
\smallskip

\begin{tabular}{@{}l@{\hspace{3mm}}l@{\hspace{3mm}}l}
 $\partial \,\Delta_{d+1}$ & -- & boundary complex of the $(d+1)$-simplex,\\[.75mm]
 $C_{d+1}(n)$       & -- & cyclic $(d+1)$-polytope with $n$ vertices,\\[.75mm]
 $C_{d+1}^{\Delta}$ & -- & $(d+1)$-dimensional cross-polytope,\\[.75mm]
 $BiC(p,q;n)$       & -- & bicyclic $4$-polytope\\[.75mm]
 $TriC(p,q,r;n)$    & -- & tricyclic $6$-polytope\\[.75mm]
 $k*k$              & -- & join product of two $k$-gons,\\[.75mm]
 $K_1*K_2$          & -- & join product of two complexes,\\[.75mm]
 $K_1\,\wr\,K_2$       & -- & wreath product of two complexes \cite{JoswigLutz2005},\\[.75mm]
 minimal            & -- & triangulation is vertex-minimal,\\[.75mm]
 tight              & -- & triangulation is tight \cite{KuehnelLutz1999},\\[.75mm]
 nncs               & -- & nearly neighborly centrally symmetric
                           sphere \cite{Lutz2004apre},\\[.75mm]
 no flip            & -- & example has no non-trivial bistellar
                           flip \cite{KuehnelLutz2003pre}.
\end{tabular}


\begin{landscape}

\small
\defaultaddspace=.1em

\setlength{\LTleft}{0pt}
\setlength{\LTright}{0pt}
\begin{longtable}{@{\extracolsep{\fill}}rcclllll@{}}
\caption{\protect\parbox[t]{15cm}{Vertex-transitive combinatorial $2$-manifolds.}}\label{tbl:2man}
\\\toprule
 \addlinespace
 \addlinespace
 \addlinespace
 $n$ & Or. &   Gen.  & $f$-vector & Group             & Type               & List of orbits & Remarks \\ \midrule
\endfirsthead
\caption{\protect\parbox[t]{15cm}{Vertex-transitive combinatorial $2$-manifolds (continued).}}
\\\toprule
 \addlinespace
 \addlinespace
 \addlinespace
 $n$ & Or. &   Gen.  & $f$-vector & Group             & Type               & List of orbits & Remarks \\ \midrule
\endhead
\bottomrule
\endfoot
 \addlinespace
 \addlinespace
 \addlinespace
  4  & $+$ &    0    & (\underline{6},\underline{4})  & $S_4$ & $\mbox{}^2\hspace{.3pt}4^{\,5}_{\,1}$   & $123_{\,4}$ & tetrahedron, regular \\
 \addlinespace
  6  & $+$ &    0    &    (12,8)  & $[2^3]S_3=2wrS_3$ & $\mbox{}^2\hspace{.3pt}6^{\,11}_{\,1}$   & $123_{\,8}$   & octahedron, regular \\
 \addlinespace
     & $-$ &    1    & (\underline{15},10) & $A_5$    & $\mbox{}^2\hspace{.3pt}6^{\,12}_{\,1}$   &  $123_{\,10}$ & ${\mathbb R}{\bf P}^{\,2}_{\,6}$, regular \\
 \addlinespace
  7  & $+$ &    1    & (\underline{21},14) & $7\!:\!6$& $\mbox{}^2\hspace{.3pt}7^{\,4}_{\,1}$    &  $124_{\,14}$ & M\"obius' torus, \cite{Csaszar1949}, \cite[$M^2$]{Kuehnel1986a-series}, \\[-.2mm]
     &     &                                 &                        &     &               && \cite{KuehnelLassmann1988-dtori}, \cite[$M^2_1$]{KuehnelLassmann1996-bundle}, \cite{Moebius1886}\\
 \addlinespace
  8  & $+$ &    1    & (24,16)    & $t8n15(32)$       & $\mbox{}^2\hspace{.3pt}8^{\,15}_{\,1}$   &  $123_{\,16}$ & \cite[$M^2_1(8)$]{KuehnelLassmann1996-bundle} \\
 \addlinespace
  9  & $+$ &    1    & (27,18)    & $D_9$             & $\mbox{}^2\hspace{.3pt}9^{\,3}_{\,1}$    &  $124_{\,18}$ & \cite[$M^2_1(9)$]{KuehnelLassmann1996-bundle} \\
 \addlinespace
     &     &         &            & ${\mathbb Z}_3^{\,2}\!:\!D_6$ & $\mbox{}^2\hspace{.3pt}9^{\,18}_{\,1}$   & $136_{\,18}$ & regular, \cite{Altshuler1973}, \cite{Brahana1926}, \cite[Ch.~8]{CoxeterMoser1957}, \cite{Wilson1976}  \\
\end{longtable}

\end{landscape}


\begin{landscape}

\small
\defaultaddspace=.1em

\setlength{\LTleft}{0pt}
\setlength{\LTright}{0pt}
\begin{longtable}{@{\extracolsep{\fill}}rllllll@{}}
\caption{\protect\parbox[t]{15cm}{Vertex-transitive combinatorial $3$-manifolds.}}\label{tbl:3man}
\\\toprule
 \addlinespace
 \addlinespace
 \addlinespace
 $n$ & Manifold & $f$-vector     & Group & Type & List of orbits & Remarks \\ \midrule
\endfirsthead
\caption{\protect\parbox[t]{15cm}{Vertex-transitive combinatorial $3$-manifolds (continued).}}
\\\toprule
 \addlinespace
 \addlinespace
 \addlinespace
 $n$ & Manifold & $f$-vector     & Group & Type & List of orbits & Remarks \\ \midrule
\endhead
\bottomrule
\endfoot
 \addlinespace
 \addlinespace
 \addlinespace
  5  & $S^3$              & (\underline{10},\underline{10},\underline{5}) & $S_5$ & $\mbox{}^3\hspace{.3pt}5^{\,5}_{\,1}$ & $1234_{\,5}$ & $\partial \,\Delta_4$, regular \\
 \addlinespace
  6  & $S^3$              & (\underline{15},18,9) & $S_3wr2$ & $\mbox{}^3\hspace{.3pt}6^{\,13}_{\,1}$ & $1234_{\,9}$ & $\partial C_4(6)=3*3$, \cite[${\rm I}_6$]{KuehnelLassmann1985-di} \\
 \addlinespace
  7  & $S^3$              & (\underline{21},28,14) & $D_7$      & $\mbox{}^3\hspace{.3pt}7^{\,2}_{\,1}$ & $1234_{\,7}$\, $1245_{\,7}$ & $\partial C_4(7)$, \cite[${\rm I}_7$]{KuehnelLassmann1985-di}\\
 \addlinespace
  8  & $S^3$              & (24,32,16) & $2wrS_4$ & $\mbox{}^3\hspace{.3pt}8^{\,44}_{\,1}$ & $1234_{\,16}$ & $\partial C_4^{\Delta}=\partial BiC(1,3;8)$ \\
     &                    &            &          &               && $=4*4$, regular, nncs\\
 \addlinespace
     &                    & (\underline{28},40,20) & $D_8$ & $\mbox{}^3\hspace{.3pt}8^{\,6}_{\,1}$ & $1234_{\,8}$\, $1245_{\,8}$\, $1256_{\,4}$ & $\partial C_4(8)$, \cite[${\rm I}_8$]{KuehnelLassmann1985-di} \\
 \addlinespace
  9  & $S^3$              & (\underline{36},54,27) & $D_9$ & $\mbox{}^3\hspace{.3pt}9^{\,3}_{\,1}$ & $1234_{\,9}$\, $1245_{\,9}$\, $1256_{\,9}$ & $\partial C_4(9)$, \cite[${\rm I}_9$]{KuehnelLassmann1985-di}\\
 \addlinespace
     & $S^2\hbox{$\times\hspace{-1.62ex}\_\hspace{-.4ex}\_\hspace{.7ex}$}S^1$ & (\underline{36},54,27) & $D_9$      & $\mbox{}^3\hspace{.3pt}9^{\,3}_{\,2}$ & $1235_{\,18}$\, $1245_{\,9}$ & minimal, tight, \\
\end{longtable}

\end{landscape}


\begin{landscape}

\small
\defaultaddspace=.1em

\setlength{\LTleft}{0pt}
\setlength{\LTright}{0pt}
\begin{longtable}{@{\extracolsep{\fill}}rllllll@{}}
\caption{\protect\parbox[t]{15cm}{Vertex-transitive combinatorial $4$-manifolds.}}\label{tbl:4man}
\\\toprule
 \addlinespace
 \addlinespace
 \addlinespace
 $n$ & Manifold & $f$-vector     & Group & Type & List of orbits & Remarks \\ \midrule
\endfirsthead
\caption{\protect\parbox[t]{15cm}{Vertex-transitive combinatorial $4$-manifolds (continued).}}
\\\toprule
 \addlinespace
 \addlinespace
 \addlinespace
 $n$ & Manifold & $f$-vector     & Group & Type & List of orbits & Remarks \\ \midrule
\endhead
\bottomrule
\endfoot
 \addlinespace
 \addlinespace
 \addlinespace
  6  & $S^4$              & (\underline{15},\underline{20},\underline{15},\underline{6}) & $S_6$ & $\mbox{}^4\hspace{.3pt}6^{\,16}_{\,1}$ & $12345_{\,6}$ & $\partial \,\Delta_5$, regular \\
 \addlinespace
  9  & ${\mathbb C}{\bf P}^{\,2}$ &(\underline{36},\underline{84},90,36)& ${\mathbb Z}_3^{\,2}\!:\!{\mathbb Z}_6$   
                                                     & $\mbox{}^4\hspace{.3pt}9^{\,13}_{\,1}$   & $12345_{\,9}$\, $12347_{\,27}$ & ${\mathbb C}{\bf P}^{\,2}_{\,9}$, minimal, tight, \\
     &                    &               &          &                    && \cite{Kuehnel1995-book}, \cite{KuehnelBanchoff1983}, \cite{KuehnelLassmann1983-unique}, \cite{KuehnelLutz1999} \\
 \addlinespace
 10  & $S^4$              & (40,80,80,32) & $2wrS_5$ & $\mbox{}^4\hspace{.3pt}10^{\,39}_{\,1}$  & $12345_{\,32}$ & $\partial C_5^{\Delta}$, regular, nncs \\
 \addlinespace
     &                    &(\underline{45},100,105,42)& ${\mathbb Z}_5$ & $\mbox{}^4\hspace{.3pt}10^{\,1}_{\,1}$   & $12345_{\,10}$\, $12356_{\,10}$\, $12367_{\,10}$\, $12379_{\,10}$\, $13579_{\,2}$ & \\
 \addlinespace
     &                    &               & $D_5$    & $\mbox{}^4\hspace{.3pt}10^{\,2}_{\,1}$   & $12345_{\,10}$\, $12356_{\,10}$\, $12367_{\,10}$\, $12379_{\,10}$\, $13579_{\,2}$ & \\
 \addlinespace
     &                    &               & $\frac{1}{2}[5\!:\!4]2$ & $\mbox{}^4\hspace{.3pt}10^{\,4}_{\,1}$   & $12345_{\,10}$\, $12347_{\,20}$\, $12359_{\,10}$\, $13579_{\,2}$ & \\
 \addlinespace
 11  & $S^3\!\times\!S^1$ &(\underline{55},110,110,44)& $D_{11}$
                                          & $\mbox{}^4\hspace{.3pt}11^{\,2}_{\,1}$   & $12346_{\,22}$\, $12356_{\,22}$ & min., tight, \cite[$M^4$]{Kuehnel1986a-series}, \\
\end{longtable}

\end{landscape}


\begin{landscape}

\small
\defaultaddspace=.1em

\setlength{\LTleft}{0pt}
\setlength{\LTright}{0pt}
\begin{longtable}{@{\extracolsep{\fill}}rllllll@{}}
\caption{\protect\parbox[t]{15cm}{Vertex-transitive combinatorial $5$-manifolds.}}\label{tbl:5man}
\\\toprule
 \addlinespace
 \addlinespace
 \addlinespace
 $n$ & Manifold & $f$-vector     & Group & Type & List of orbits & Remarks \\ \midrule
\endfirsthead
\caption{\protect\parbox[t]{15cm}{Vertex-transitive combinatorial $5$-manifolds (continued).}}
\\\toprule
 \addlinespace
 \addlinespace
 \addlinespace
 $n$ & Manifold & $f$-vector     & Group & Type & List of orbits & Remarks \\ \midrule
\endhead
\bottomrule
\endfoot
 \addlinespace
 \addlinespace
 \addlinespace
  7  & $S^5$ & (\underline{21},\underline{35},\underline{35}, & $S_7$ & $\mbox{}^5\hspace{.3pt}7^{\,7}_{\,1}$ & $123456_{\,7}$ & $\partial \,\Delta_6$, regular \\
     &       &\underline{21},\underline{7})&&&& \\
 \addlinespace
  8  & $S^5$ &(\underline{28},\underline{56},68, & $S_4wr2$ & $\mbox{}^5\hspace{.3pt}8^{\,47}_{\,1}$   & $123456_{\,16}$& $\partial C_6(8)=(\partial \,\Delta_3)^{*2}$ \\
     &       & 48,16)         &              &       &                & $=\partial\Delta_1\wr\partial C_2(4)$ \cite{JoswigLutz2005} \\
 \addlinespace
  9  & $S^5$ &(\underline{36},81,108, & $S_3wrS_3$ & $\mbox{}^5\hspace{.3pt}9^{\,31}_{\,1}$   & $123467_{\,27}$& $3*3*3$ \\
     &       & 81,27)         &              &       &                & $=\partial TriC(1,2,4;9)$ \\
 \addlinespace
     &       &(\underline{36},\underline{84},117,   & $D_9$       & $\mbox{}^5\hspace{.3pt}9^{\,3}_{\,1}$    & $123456_{\,9}$\, $123467_{\,18}$\, $124578_{\,3}$ & $\partial C_6(9)$ \\
     &       & 90,30)&&&& \\
\end{longtable}

\end{landscape}


%


\begin{landscape}

\small
\defaultaddspace=.1em

\setlength{\LTleft}{0pt}
\setlength{\LTright}{0pt}
\begin{longtable}{@{\extracolsep{\fill}}rllllll@{}}
\caption{\protect\parbox[t]{19cm}{Vertex-transitive combinatorial $6$-manifolds.}}\label{tbl:6man}
\\\toprule
 \addlinespace
 \addlinespace
 \addlinespace
 $n$ & Manifold & $f$-vector     & Group & Type & List of orbits & Remarks \\ \midrule
\endfirsthead
\caption{\protect\parbox[t]{15cm}{Vertex-transitive combinatorial $6$-manifolds (continued).}}
\\\toprule
 \addlinespace
 \addlinespace
 \addlinespace
 $n$ & Manifold & $f$-vector     & Group & Type & List of orbits & Remarks \\ \midrule
\endhead
\bottomrule
\endfoot
 \addlinespace
 \addlinespace
 \addlinespace
   8 & $S^6$ & (\underline{28},\underline{56},\underline{70},\underline{56},\underline{28},\underline{8}) & $S_8$ & $\mbox{}^6\hspace{.3pt}8^{\,50}_{\,1}$ & $1234567_{\,8}$ & $\partial \,\Delta_7$, regular \\
 \addlinespace
  10 & $S^6$ &(\underline{45},\underline{120},205,222,140,40) & $2wrD_5$ & $\mbox{}^6\hspace{.3pt}10^{\,23}_{\,1}$ & $1234567_{\,40}$ & $\partial\Delta_1\wr\partial C_2(5)$ \cite{JoswigLutz2005}\\
 \addlinespace
  14 & $S^6$ &(84,280,560,672,448,128)& $2wrS_7$   & $\mbox{}^6\hspace{.3pt}14^{\,57}_{\,1}$ & $1234567_{\,128}$ & $\partial C_7^{\Delta}$, regular, nncs \\
 \addlinespace
     &       &(\underline{91},322,665,812,546,156)& $2[\frac{1}{2}]7\!:\!6$ & $\mbox{}^6\hspace{.3pt}14^{\,4}_{\,1}$ & $1234567_{\,42}$\, $123457\,13_{\,42}$\, $123467\,12_{\,14}$ & \\ 
     &       &                &                      &                  & $12356\,11\,12_{\,42}$\, $12357\,11\,13_{\,14}$\, $13579\,11\,13_{\,2}$ & \\
 \addlinespace
 &&(\underline{91},\underline{364},875,1190,840,240)& $2[\frac{1}{2}]7\!:\!6$ & $\mbox{}^6\hspace{.3pt}14^{\,4}_{\,2}$ & $1234568_{\,42}$\, $123458\,14_{\,42}$\, $12345\,11\,13_{\,42}$ & \\
     &       &                &                      &                  & $123567\,10_{\,42}$\, $123567\,13_{\,42}$\, $12356\,11\,14_{\,14}$ & \\ 
     &       &                &                      &                  & $12357\,11\,13_{\,14}$\, $13579\,11\,13_{\,2}$ & \\
\end{longtable}

\end{landscape}


\begin{landscape}

\small
\defaultaddspace=.2em

\setlength{\LTleft}{0pt}
\setlength{\LTright}{0pt}
\begin{longtable}{@{\extracolsep{\fill}}rllllll@{}}
\caption{\protect\parbox[t]{15cm}{Vertex-transitive combinatorial $7$-manifolds.}}\label{tbl:7man}
\\\toprule
 \addlinespace
 \addlinespace
 $n$ & Manifold & $f$-vector     & Group & Type & List of orbits & Remarks \\ \midrule
\endfirsthead
\caption{\protect\parbox[t]{15cm}{Vertex-transitive combinatorial $7$-manifolds (continued).}}
\\\toprule
 \addlinespace
 \addlinespace
 $n$ & Manifold & $f$-vector     & Group & Type & List of orbits & Remarks \\ \midrule
\endhead
\bottomrule
\endfoot
 \addlinespace
 \addlinespace
 \addlinespace
  9  & $S^7$ & (\underline{36},\underline{84},\underline{126},\underline{126}, & $S_9$ & $\mbox{}^7\hspace{.3pt}9^{\,34}_{\,1}$ & $12345678_{\,9}$ & $\partial \,\Delta_8$, regular \\
     &       & \underline{84},\underline{36},\underline{9}) &&&& \\
 \addlinespace
 10  & $S^7$ &(\underline{45},\underline{120},\underline{210},250, & $S_5wr2$ & $\mbox{}^7\hspace{.3pt}10^{\,43}_{\,1}$  & $12345678_{\,25}$ & $\partial C_8(10)=(\partial \,\Delta_4)^{*2}$ \\
     &       & 200,100,25)    &                      &        &                & \\
 \addlinespace
 11  & $S^7$ &(\underline{55},\underline{165},\underline{330},451, & $D_{11}$ & $\mbox{}^7\hspace{.3pt}11^{\,2}_{\,1}$   & $12345678_{\,11}$\, $12345689_{\,22}$\, $12346789_{\,11}$ & $\partial C_8(11)$ \\
     &       & 407,220,55)    &                      &                  & $1234679\,10_{\,11}$ & \\
 \addlinespace
 12  & $S^7$ &(\underline{66},216,459,648,& $S_3wrS_4$              & $\mbox{}^7\hspace{.3pt}12^{\,289}_{\,1}$ & $12345678_{\,81}$ & $3*3*3*3$ \\
     &       & 594,324,81) &&&& \\
\end{longtable}

\end{landscape}


\begin{landscape}

\small
\defaultaddspace=.2em

\setlength{\LTleft}{0pt}
\setlength{\LTright}{0pt}
\begin{longtable}{@{\extracolsep{\fill}}rllllll@{}}
\caption{\protect\parbox[t]{15cm}{Vertex-transitive combinatorial $8$-manifolds.}}\label{tbl:8man}
\\\toprule
 \addlinespace
 \addlinespace
 $n$ & Manifold & $f$-vector     & Group & Type & List of orbits & Remarks \\ \midrule
\endfirsthead
\caption{\protect\parbox[t]{15cm}{Vertex-transitive combinatorial $8$-manifolds (continued).}}
\\\toprule
 \addlinespace
 \addlinespace
 $n$ & Manifold & $f$-vector     & Group & Type & List of orbits & Remarks \\ \midrule
\endhead
\bottomrule
\endfoot
 \addlinespace
 \addlinespace
 \addlinespace
 10  & $S^8$ & (\underline{45},\underline{120},\underline{210},\underline{252}, & $S_{10}$ & $\mbox{}^8\hspace{.3pt}10^{\,45}_{\,1}$ & $123456789_{\,10}$ & $\partial \,\Delta_9$, regular \\
     &       & \underline{210},\underline{120},\underline{45},\underline{10}) &&&& \\
 \addlinespace
 12  & $S^8$ &(\underline{66},\underline{220},492,768, &$S_4wrS_3$& $\mbox{}^8\hspace{.3pt}12^{\,294}_{\,1}$ & $123456789_{\,64}$ & $(\partial \,\Delta_3)^{*3}$ \\
     &       & 840,624,288,64) &&&& $=\partial\Delta_1\wr\partial C_3^{\Delta}$ \cite{JoswigLutz2005}\\
 \addlinespace
 14  & $S^8$ &(\underline{91},\underline{364},987,1862, & ${\mathbb Z}_{14}$ & $\mbox{}^8\hspace{.3pt}14^{\,1}_{\,5}$ & $123456789_{\,14}$\, $12345679\,11_{\,14}$\, $1234567\,11\,13_{\,14}$ & \\ 
     &       & 2408,2032,1008,224) &                      &                  & $12345689\,10_{\,14}$\, $1234569\,10\,11_{\,14}$\, $123456\,10\,11\,12_{\,14}$ & \\ 
     &       &                &                      &                  & $123456\,11\,12\,13_{\,14}$\, $1234579\,11\,13_{\,14}$\, $12346789\,11_{\,14}$ & \\ 
     &       &                &                      &                  & $1234678\,11\,13_{\,14}$\, $1234689\,10\,11_{\,14}$\, $123468\,10\,11\,13_{\,14}$ & \\
\end{longtable}

\end{landscape}


\begin{landscape}

\small
\defaultaddspace=.2em

\setlength{\LTleft}{0pt}
\setlength{\LTright}{0pt}
\begin{longtable}{@{\extracolsep{\fill}}rllllll@{}}
\caption{\protect\parbox[t]{15cm}{Vertex-transitive combinatorial $9$-manifolds.}}\label{tbl:9man}
\\\toprule
 \addlinespace
 \addlinespace
 $n$ & Manifold & $f$-vector     & Group & Type & List of orbits & Remarks \\ \midrule
\endfirsthead
\caption{\protect\parbox[t]{15cm}{Vertex-transitive combinatorial $9$-manifolds (continued).}}
\\\toprule
 \addlinespace
 \addlinespace
 $n$ & Manifold & $f$-vector     & Group & Type & List of orbits & Remarks \\ \midrule
\endhead
\bottomrule
\endfoot
 \addlinespace
 \addlinespace
 \addlinespace
 11  & $S^9$ & (\underline{55},\underline{165},\underline{330},\underline{462},\underline{462}, & $S_{11}$ & $\mbox{}^9\hspace{.3pt}11^{\,8}_{\,1}$ & $123456789\,10_{\,11}$ & $\partial \,\Delta_{10}$, regular \\
     &       & \underline{330},\underline{165},\underline{55},\underline{11}) &&&& \\
 \addlinespace
 12  & $S^9$ &(\underline{66},\underline{220},\underline{495},\underline{792},922,& $S_6wr2$ & $\mbox{}^9\hspace{.3pt}12^{\,299}_{\,1}$ & $123456789\,10_{\,36}$ & $\partial C_{10}(12)=(\partial \,\Delta_5)^{*2}$ \\
     &       &                               &              &       &                & $=\partial\Delta_1\wr\partial C_4(6)$ \\
     &       &                               &              &       &                & $=\partial\Delta_2\wr\partial C_2(4)$ \cite{JoswigLutz2005} \\
 \addlinespace
 13  & $S^9$ &(\underline{78},\underline{286},\underline{715},\underline{1287},1703, &    $D_{13}$      & $\mbox{}^9\hspace{.3pt}13^{\,2}_{\,1}$ & $123456789\,10_{\,13}$\, $12345678\,10\,11_{\,26}$  & $\partial C_{10}(13)$ \\
     &       & 1638,1092,455,91)               &                         &                    & $12345689\,10\,11_{\,26}$\, $12345689\,11\,12_{\,13}$ & \\
     &       &                &                         &                    & $12346789\,11\,12_{\,13}$ & \\
\end{longtable}

\end{landscape}


\begin{landscape}

\small
\defaultaddspace=.15em

\setlength{\LTleft}{0pt}
\setlength{\LTright}{0pt}
\begin{longtable}{@{\extracolsep{\fill}}rllllll@{}}
\caption{\protect\parbox[t]{19cm}{Vertex-transitive combinatorial $10$-manifolds.}}\label{tbl:10man}
\\\toprule
 \addlinespace
 \addlinespace
 $n$ & Manifold & $f$-vector     & Group & Type & List of orbits & Remarks \\ \midrule
\endhead
\bottomrule
\endfoot
 \addlinespace
 \addlinespace
 \addlinespace
 12  & $S^{10}$ & (\underline{66},\underline{220},\underline{495},\underline{792},\underline{924},\underline{792},\underline{495},\underline{220},\underline{66},\underline{12}) & $S_{12}$ & $\mbox{}^{10}\hspace{.3pt}12^{\,301}_{\,1}$ & $123456789\,10\,11_{\,12}$ & $\partial \,\Delta_{11}$, regular \\ 
 \addlinespace
 14  & $S^{10}$ & (\underline{91},\underline{364},\underline{1001},\underline{2002},2996,3376,2814,1652,616,112) & $2wrD_7$ & $\mbox{}^{10}\hspace{.3pt}14^{\,38}_{\,1}$ & $123456789\,10\,11_{\,56}$ & $\partial\Delta_1\wr\partial C_4(7)$ \cite{JoswigLutz2005}\\
     &          &                &                         &                    & $123456789\,11\,12_{\,56}$ & \\
 \addlinespace
 \addlinespace
 \addlinespace
\end{longtable}

\bigskip


\small
\defaultaddspace=.15em

\setlength{\LTleft}{0pt}
\setlength{\LTright}{0pt}
\begin{longtable}{@{\extracolsep{\fill}}rllllll@{}}
\caption{\protect\parbox[t]{19cm}{Vertex-transitive combinatorial $11$-manifolds.}}\label{tbl:11man}
\\\toprule
 \addlinespace
 \addlinespace
 $n$ & Manifold & $f$-vector     & Group & Type & List of orbits & Remarks \\ \midrule
\endhead
\bottomrule
\endfoot
 \addlinespace
 \addlinespace
 \addlinespace
 13  & $S^{11}$ & (\underline{78},\underline{286},\underline{715},\underline{1287},\underline{1716},\underline{1716}, & $S_{13}$ & $\mbox{}^{11}\hspace{.3pt}13^{\,9}_{\,1}$ & $123456789\,10\,11\,12_{\,13}$ & $\partial \,\Delta_{12}$, regular \\
     &          & \underline{1287},\underline{715},\underline{286},\underline{78},\underline{13}) &&&& \\ 
 \addlinespace
 14  & $S^{11}$ &(\underline{91},\underline{364},\underline{1001},\underline{2002},\underline{3003},3430,& $S_7wr2$ & $\mbox{}^{11}\hspace{.3pt}14^{\,61}_{\,1}$ & $123456789\,10\,11\,12_{\,49}$ & $\partial C_{12}(14)=(\partial \,\Delta_6)^{*2}$ \\
     &          & 2989,1960,931,294,49)               &              &       &                &  \\
 \addlinespace
 15  & $S^{11}$ &(\underline{105},\underline{455},\underline{1365},3000,4975,6300,& $S_5wrS_3$ & $\mbox{}^{11}\hspace{.3pt}15^{\,102}_{\,1}$ & $123456789\,10\,11\,12_{\,125}$ & $(\partial \,\Delta_4)^{*3}$ \\
     &          & 6075,4375,2250,750,125) &&&& \\
 \addlinespace
     &          &(\underline{105},\underline{455},\underline{1365},\underline{3003},5000,6390,& $S_3wrD_5$ & $\mbox{}^{11}\hspace{.3pt}15^{\,86}_{\,1}$ & $123456789\,10\,11\,12_{\,135}$ & $\partial\Delta_2\wr\partial C_2(5)$ \cite{JoswigLutz2005}\\
     &          & 6255,4590,2403,810,135) &&&& \\
\end{longtable}

\end{landscape}

\begin{landscape}

\small
\defaultaddspace=.2em

\setlength{\LTleft}{0pt}
\setlength{\LTright}{0pt}
\begin{longtable}{@{\extracolsep{\fill}}lll@{}}
\caption{\protect\parbox[t]{15cm}{List of generators for the group actions in Table~\ref{tbl:2man} -- Table~\ref{tbl:11man}.}}\label{tbl:generators}
\\\toprule
 \addlinespace
 \addlinespace
     Action   &        Group        & Generators \\ \midrule
\endfirsthead
\caption{\protect\parbox[t]{15cm}{List of generators (continued).}}
\\\toprule
 \addlinespace
 \addlinespace
     Action   &        Group        & Generators \\ \midrule
\endhead
\bottomrule
\endfoot
 \addlinespace
 \addlinespace
 \addlinespace
   $6^{\,11}$ &$[2^3]S_3\!=\!2wrS_3$&  (3,6), (1,3,5)(2,4,6), (1,5)(2,4) \\
 \addlinespace
   $6^{\,12}$ & $A_5$               &  (1,2,3,4,6), (1,4)(5,6)\\
 \addlinespace
   $6^{\,13}$ & $[S_3^{\,2}]2\!=\!S_3wr2$  &  (2,4,6), (2,4), (1,4)(2,5)(3,6) \\
 \addlinespace
    $7^{\,4}$ & $7\!:\!6$           &  (1,2,3,4,5,6,7), (1,3,2,6,4,5)\\
 \addlinespace
   $8^{\,15}$ & $t8n15(32)$         &  (1,2,3,4,5,6,7,8), (1,5)(3,7), (1,6)(2,5)(3,4)(7,8)\\
 \addlinespace
   $8^{\,44}$ &$[2^4]S_4\!=\!2wrS_4$&  (4,8), (1,8)(4,5), (1,2,3,8)(4,5,6,7)\\
 \addlinespace
   $8^{\,47}$ & $[S_4^{\,2}]2\!=\!S_4wr2$  &  (1,2,3,8), (2,3), (1,5)(2,6)(3,7)(4,8) \\
 \addlinespace
    $9^{\,4}$ & $S_3\times{\mathbb Z}_3$ & (1,2,9)(3,4,5)(6,7,8), (1,2)(4,5)(7,8), (1,4,7)(2,5,8)(3,6,9) \\
\end{longtable}

\end{landscape}


\pagebreak

\bibliography{}

\begin{thebibliography}{10}

\bibitem{Altshuler1973}
A.~Altshuler.
\newblock Construction and enumeration of regular maps on the torus.
\newblock {\em Discrete Math.} {\bf 4}, 201--217 (1973).

\bibitem{Altshuler1974}
A.~Altshuler.
\newblock Combinatorial $3$-manifolds with few vertices.
\newblock {\em J.\ Comb.\ Theory, Ser.\ A\/} {\bf 16}, 165--173 (1974).

\bibitem{Altshuler1976}
A.~Altshuler.
\newblock A peculiar triangulation of the $3$-sphere.
\newblock {\em Proc.\ Am.\ Math.\ Soc.} {\bf 54}, 449--452 (1976).

\bibitem{Altshuler1977}
A.~Altshuler.
\newblock Neighborly $4$-polytopes and neighborly combinatorial
  $3$-man\-i\-folds with ten vertices.
\newblock {\em Can.\ J.\ Math.} {\bf 29}, 400--420 (1977).

\bibitem{AltshulerBokowskiSchuchert1996}
A.~Altshuler, J.~Bokowski, and P.~Schuchert.
\newblock Neighborly $2$-manifolds with $12$~vertices.
\newblock {\em J.\ Comb.\ Theory, Ser.\ A\/} {\bf 75}, 148--162 (1996).

\bibitem{AltshulerBrehm1992}
A.~Altshuler and U.~Brehm.
\newblock Neighborly maps with few vertices.
\newblock {\em Discrete Comput.\ Geom.} {\bf 8}, 93--104 (1992).

\bibitem{AltshulerSteinberg1973}
A.~Altshuler and L.~Steinberg.
\newblock Neighborly $4$-polytopes with $9$~vertices.
\newblock {\em J.\ Comb.\ Theory, Ser.\ A\/} {\bf 15}, 270--287 (1973).

\bibitem{AltshulerSteinberg1974}
A.~Altshuler and L.~Steinberg.
\newblock Neighborly combinatorial $3$-manifolds with $9$~ver\-tices.
\newblock {\em Discrete Math.} {\bf 8}, 113--137 (1974).

\bibitem{AltshulerSteinberg1976}
A.~Altshuler and L.~Steinberg.
\newblock An enumeration of combinatorial $3$-mani\-folds with nine vertices.
\newblock {\em Discrete Math.} {\bf 16}, 91--108 (1976).

\bibitem{Barden1965}
D.~Barden.
\newblock Simply connected five-manifolds.
\newblock {\em Ann.\ Math.} {\bf 82}, 365--385 (1965).

\bibitem{BjoernerLutz2000}
A.~Bj\"orner and F.~H. Lutz.
\newblock Simplicial manifolds, bistellar flips and a $16$-ver\-tex
  triangulation of the {Poincar\'{e}} homology $3$-sphere.
\newblock {\em Exp.\ Math.} {\bf 9}, 275--289 (2000).

\bibitem{Bokowski1989}
J.~Bokowski.
\newblock A geometric realization without self-intersections does exist for
  {Dyck's} regular map.
\newblock {\em Discrete Comput.\ Geom.} {\bf 4}, 583--589 (1989).

\bibitem{BokowskiGarms1987}
J.~Bokowski and K.~Garms.
\newblock Altshuler's sphere {$M^{10}_{425}$} is not polytopal.
\newblock {\em Eur.\ J.\ Comb.} {\bf 8}, 227--229 (1987).

\bibitem{BokowskiWills1988}
J.~Bokowski and J.~M. Wills.
\newblock Regular polyhedra with hidden symmetries.
\newblock {\em Math.\ Intell.} {\bf 10}, No. 1, 27--32 (1988).

\bibitem{Bose1939}
R.~C. Bose.
\newblock On the construction of balanced incomplete block designs.
\newblock {\em Ann.\ Eugenics\/} {\bf 9}, 353--399 (1939).

\bibitem{Brahana1926}
H.~R. Brahana.
\newblock Regular maps on an anchor ring.
\newblock {\em Am.\ J.\ Math.} {\bf 48}, 225--240 (1926).

\bibitem{Brehm-pers}
U.~Brehm.
\newblock Personal communication, 1999.

\bibitem{Brehm1987a}
U.~Brehm.
\newblock Maximally symmetric polyhedral realizations of {D}yck's regular map.
\newblock {\em Mathematika\/} {\bf 34}, 229--236 (1987).

\bibitem{BrehmKuehnel1986}
U.~Brehm and W.~K\"uhnel.
\newblock A polyhedral model for {Cartan's} hypersurface in~{$S^4$}.
\newblock {\em Mathematika\/} {\bf 33}, 55--61 (1986).

\bibitem{BrehmKuehnel1987}
U.~Brehm and W.~K\"uhnel.
\newblock Combinatorial manifolds with few vertices.
\newblock {\em Topo\-logy\/} {\bf 26}, 465--473 (1987).

\bibitem{BrehmKuehnel1992}
U.~Brehm and W.~K\"uhnel.
\newblock $15$-vertex triangulations of an $8$-manifold.
\newblock {\em Math.\ Ann.} {\bf 294}, 167--193 (1992).

\bibitem{BrehmLutz2002pre}
U.~Brehm and F.~H. Lutz.
\newblock Triangulations of {S}eifert manifolds.
\newblock In preparation.

\bibitem{Butler1993}
G.~Butler.
\newblock The transitive groups of degree fourteen and fifteen.
\newblock {\em J.\ Symb.\ Comput.} {\bf 16}, 413--422 (1993).

\bibitem{ButlerMcKay1983}
G.~Butler and J.~McKay.
\newblock The transitive groups of degree up to eleven.
\newblock {\em Commun.\ Algebra\/} {\bf 11}, 863--911 (1983).

\bibitem{CasellaKuehnel2001}
M.~Casella and W.~K\"uhnel.
\newblock A triangulated {K3} surface with the minimum number of vertices.
\newblock {\em Topology\/} {\bf 40}, 753--772 (2001).

\bibitem{ConwayHulpkeMcKay1998}
J.~H. Conway, A.~Hulpke, and J.~McKay.
\newblock On transitive permutation groups.
\newblock {\em LMS J.\ Comput.\ Math.} {\bf 1}, 1--8 (1998).

\bibitem{Coxeter1943}
H.~S.~M. Coxeter.
\newblock The map-coloring of unorientable surfaces.
\newblock {\em Duke Math.\ J.} {\bf 10}, 293--304 (1943).

\bibitem{CoxeterMoser1957}
H.~S.~M. Coxeter and W.~O.~J. Moser.
\newblock {\em Generators and Relations for Discrete Groups}. Ergebnisse der
  Mathematik und ihrer Grenzgebiete~{\bf 14}.
\newblock Springer-Verlag, Berlin, 1957.
\newblock Fourth edition, 1980.

\bibitem{Csaszar1949}
A.~{Cs\'asz\'ar}.
\newblock A polyhedron without diagonals.
\newblock {\em Acta Sci.\ Math., Szeged\/} {\bf 13}, 140--142 (1949--1950).

\bibitem{DumasHeckenbachSaundersWelker2003a}
J.-G. Dumas, F.~Heckenbach, B.~D. Saunders, and V.~Welker.
\newblock Computing simplicial homology based on efficient {S}mith normal form
  algorithms.
\newblock {\em Algebra, Geometry, and Software Systems} (M.~Joswig and
  N.~Takayama, eds.),  177--206. Springer-Verlag, Berlin, 2003.

\bibitem{DumasHeckenbachSaundersWelker2003b}
J.-G. Dumas, F.~Heckenbach, B.~D. Saunders, and V.~Welker.
\newblock \emph{Simplicial {H}omology}, a (proposed) {GAP} share package,
  {V}ersion 1.4.1.
\newblock \url{http://www.cis.udel.edu/~dumas/Homology/}, 2003.

\bibitem{Dyck1880b}
W.~Dyck.
\newblock Notiz \"uber eine regul\"are {R}iemann'sche {F}l\"ache vom
  {G}eschlechte drei und die zugeh\"orige ,,{N}ormalcurve`` vierter {O}rdnung.
\newblock {\em Math.\ Ann.} {\bf 17}, 510--516 (1880).

\bibitem{Dyck1880a}
W.~Dyck.
\newblock Ueber {A}ufstellung und {U}ntersuchung von {G}ruppe und
  {I}rrationalit{\"a}t regul\"arer {R}iemann'scher {F}l\"achen.
\newblock {\em Math.\ Ann.} {\bf 17}, 473--509 (1880).

\bibitem{EellsKuiper1962b}
J.~{Eells, Jr.} and N.~H. Kuiper.
\newblock Manifolds which are like projective planes.
\newblock {\em Publ.\ Math., Inst.\ Hautes \'Etud.\ Sci.} {\bf 14}, 181--222
  (1962).

\bibitem{polymake}
E.~Gawrilow and M.~Joswig.
\newblock \texttt{polymake}: a {F}ramework for {A}nalyzing {C}onvex {P}olytopes
  and {S}implicial {C}omplexes.
\newblock \url{http://www.math.tu-berlin.de/polymake}, 1997--2004.
\newblock Version 2.1.0, with contributions by Thilo Schr\"oder and Nikolaus
  Witte.

\bibitem{Gorbatsevich1980}
V.~V. Gorbatsevich.
\newblock On compact homogeneous manifolds of low dimension.
\newblock {\em Geom.\ Metody Zadachakh Algebry Anal.} {\bf 2}, 37--60 (1980).

\bibitem{GAP4}
The~GAP Group.
\newblock {GAP} -- {G}roups, {A}lgorithms, and {P}rogramming, {V}ersion~4.4.
\newblock \url{http://www.gap-system.org}, 2004.

\bibitem{Gruenbaum1967}
B.~Gr\"unbaum.
\newblock {\em Convex Polytopes}. Pure and Applied Mathematics~{\bf 16}.
\newblock Interscience Publishers, London, 1967.
\newblock Second edition (V.~Kaibel, V.~Klee, and G.~M.~Ziegler, eds.),
  Graduate Texts in Mathematics \textbf{221}. Springer-Verlag, New York, NY,
  2003.

\bibitem{Heckenbach1997}
F.~Heckenbach.
\newblock \emph{{D}ie {M}\"obiusfunktion und {H}omologien auf partiell
  ge\-ord\-neten {M}engen}.
\newblock Thesis for Diploma at University Erlangen-Nuremberg, 1997. Computer
  program \texttt{homology}, \url{http://www.mi.uni-erlangen.de/~heckenb/}.

\bibitem{Hulpkepre}
A.~Hulpke.
\newblock Constructing transitive permutation groups.
\newblock \emph{J.\ Symbolic Comput.}, to appear.

\bibitem{Hulpke1996}
A.~Hulpke.
\newblock {\em Konstruktion transitiver Permutationsgruppen. {\rm
  Dissertation}}.
\newblock Verlag der Augustinus-Buchhandlung, Aachen, 1996.

\bibitem{JoswigLutz2005}
M.~Joswig and F.~H. Lutz.
\newblock One-point suspensions and wreath products of polytopes and spheres.
\newblock {\em J.\ Comb.\ Theory, Ser.\ A\/} {\bf 110}, 193--216 (2005).

\bibitem{KimmerleKouzoudi2003}
W.~Kimmerle and E.~Kouzoudi.
\newblock Doubly transitive automorphism groups of combinatorial surfaces.
\newblock {\em Discrete Comput.\ Geom.} {\bf 29}, 445--457 (2003).

\bibitem{Klaus1988}
S.~Klaus.
\newblock \emph{Einfach-zusammenh\"angende kompakte homogene {R}\"aume bis zur
  {D}imension neun}.\, {D}iplomarbeit.
\newblock Johannes Gutenberg Universit\"at Mainz, 1988, 98 pages.

\bibitem{Kreck2001}
M.~Kreck.
\newblock An inverse to the {P}oincar\'e conjecture.
\newblock {\em Arch.\ Math.} {\bf 77}, 98--106 (2001).

\bibitem{Kuehnel1986a-series}
W.~K\"uhnel.
\newblock Higherdimensional analogues of {C}s\'asz\'ar's torus.
\newblock {\em Result.\ Math.} {\bf 9}, 95--106 (1986).

\bibitem{Kuehnel1990-few}
W.~K\"uhnel.
\newblock Triangulations of manifolds with few vertices.
\newblock {\em Advances in Differential Geometry and Topology} (F.~Tricerri,
  ed.),  59--114. World Scientific, Singapore, 1990.

\bibitem{Kuehnel1995-book}
W.~K\"uhnel.
\newblock {\em Tight Polyhedral Submanifolds and Tight Triangulations}. Lecture
  Notes in Mathematics~{\bf 1612}.
\newblock Springer-Verlag, Berlin, 1995.

\bibitem{KuehnelBanchoff1983}
W.~K\"uhnel and T.~F. Banchoff.
\newblock The $9$-vertex complex projective plane.
\newblock {\em Math.\ Intell.} {\bf 5}, No. 3, 11--22 (1983).

\bibitem{KuehnelLassmann1983-unique}
W.~K\"uhnel and G.~Lassmann.
\newblock The unique $3$-neighborly $4$-manifold with few vertices.
\newblock {\em J.\ Comb.\ Theory, Ser.\ A\/} {\bf 35}, 173--184 (1983).

\bibitem{KuehnelLassmann1984-3torus}
W.~K\"uhnel and G.~Lassmann.
\newblock The rhombidodecahedral tessellation of $3$-space and a particular
  $15$-vertex triangulation of the $3$-dimensional torus.
\newblock {\em Manuscr.\ Math.} {\bf 49}, 61--77 (1984).

\bibitem{KuehnelLassmann1985-di}
W.~K\"uhnel and G.~Lassmann.
\newblock Neighborly combinatorial $3$-manifolds with dihedral automorphism
  group.
\newblock {\em Isr.\ J.\ Math.} {\bf 52}, 147--166 (1985).

\bibitem{KuehnelLassmann1988-dtori}
W.~K\"uhnel and G.~Lassmann.
\newblock Combinatorial $d$-tori with a large symmetry group.
\newblock {\em Discrete Comput.\ Geom.} {\bf 3}, 169--176 (1988).

\bibitem{KuehnelLassmann1996-bundle}
W.~K\"uhnel and G.~Lassmann.
\newblock Permuted difference cycles and triangulated sphere bundles.
\newblock {\em Discrete Math.} {\bf 162}, 215--227 (1996).

\bibitem{KuehnelLutz1999}
W.~K\"uhnel and F.~H. Lutz.
\newblock A census of tight triangulations.
\newblock {\em Periodica Math.\ Hung.} {\bf 39}, 161--183 (1999).

\bibitem{KuehnelLutz2003pre}
W.~K\"uhnel and F.~H. Lutz\mbox{}\mbox{}.
\newblock Triangulated {M}anifolds with {F}ew {V}ertices: {R}ecognition of
  {M}anifolds.
\newblock In preparation.

\bibitem{LassmannSparla2000}
G.~Lassmann and E.~Sparla.
\newblock A classification of centrally-symmetric and cyclic $12$-vertex
  triangulations of {$S^2\!\times\!S^2$}.
\newblock {\em Discrete Math.} {\bf 223}, 175--187 (2000).

\bibitem{BokowskiBremnerLutzMartin2003pre}
F.~H. Lutz.
\newblock Combinatorial $3$-manifolds with $10$ vertices.
\newblock In preparation.

\bibitem{Lutz1999}
F.~H. Lutz.
\newblock {\em Triangulated Manifolds with Few Vertices and
  Vertex-Tran\-si\-tive Group Actions. {\rm Dissertation}}.
\newblock Shaker Verlag, Aachen, 1999, 146 pages.

\bibitem{Lutz_BISTELLAR_EQUIVALENT}
F.~H. Lutz\mbox{}.
\newblock \texttt{BISTELLAR\_EQUIVALENT}, {V}ersion {F}eb/1999.
\newblock
  \url{http://www.math.tu-berlin.de/diskregeom/stellar/BISTELLAR_EQUIVALENT},
  1999.

\bibitem{Lutz_MANIFOLD_VT}
F.~H. Lutz\mbox{}.
\newblock \texttt{MANIFOLD\_VT}, {V}ersion {A}pr/2002.
\newblock \url{http://www.math.tu-berlin.de/diskregeom/stellar/MANIFOLD_VT},
  2002.

\bibitem{Lutz_BISTELLAR}
F.~H. Lutz\mbox{}.
\newblock \texttt{BISTELLAR}, {V}ersion {N}ov/2003.
\newblock \url{http://www.math.tu-berlin.de/diskregeom/stellar/BISTELLAR},
  2003.

\bibitem{Lutz_SEIFERT}
F.~H. Lutz\mbox{}.
\newblock \texttt{SEIFERT}, {V}ersion {N}ov/2003.
\newblock \url{http://www.math.tu-berlin.de/diskregeom/stellar/SEIFERT}, 2003.

\bibitem{Lutz2005apre}
F.~H. Lutz\mbox{}\mbox{}.
\newblock Enumeration and random realization of triangulated surfaces.
\newblock \url{arXiv:math.CO/0506316}, 2005, 15 pages.

\bibitem{Lutz2004apre}
F.~H. Lutz\mbox{}\mbox{}.
\newblock Triangulated {M}anifolds with {F}ew {V}ertices: {C}entrally
  {S}ymmetric {S}pheres and {P}roducts of {S}pheres.
\newblock \url{arXiv:math.MG/0404465}, 2004, 26 pages.

\bibitem{Lutz2005bpre}
F.~H. Lutz\mbox{}\mbox{}.
\newblock Triangulated {M}anifolds with {F}ew {V}ertices: {C}ombinatorial
  {M}anifolds.
\newblock \url{arXiv:math.CO/0506372}, 2005, 37 pages.

\bibitem{Lutz2004dpre}
F.~H. Lutz\mbox{}\mbox{}.
\newblock Triangulated {M}anifolds with {F}ew {V}ertices: {C}ombinatorial
  {P}seudomanifolds.
\newblock In preparation.

\bibitem{Lutz2003bpre}
F.~H. Lutz\mbox{}\mbox{}.
\newblock Triangulated {M}anifolds with {F}ew {V}ertices: {G}eometric
  $3$-{M}anifolds.
\newblock \url{arXiv:math.GT/0311116}, 2003, 48 pages.

\bibitem{Lutz2004bpre}
F.~H. Lutz\mbox{}\mbox{}.
\newblock Triangulated {M}anifolds with {F}ew {V}ertices: {V}ertex-{T}ransitive
  {T}riangulations {II}.
\newblock In preparation.

\bibitem{Matousek2003}
J.~Matou\v{s}ek.
\newblock {\em Using the Borsuk-Ulam Theorem. Lectures on Topological Methods
  in Combinatorics and Geometry}. Universitext.
\newblock Springer-Verlag, Berlin, 2003.

\bibitem{nauty}
B.~D. McKay.
\newblock \texttt{nauty}, {V}ersion 2.2.
\newblock \url{http://cs.anu.edu.au/people/bdm/nauty/}, 1994--2003.

\bibitem{Miller1896-12}
G.~A. Miller.
\newblock List of transitive substitution groups of degree twelve.
\newblock {\em Quart.\ J.\ Pure Appl.\ Math.} {\bf 28}, 193--231 (1896).
\newblock Errata: \emph{Quart.\ J.\ Pure Appl.\ Math.} \textbf{29}, 249 (1897).

\bibitem{Miller1897-13}
G.~A. Miller.
\newblock On the transitive substitution groups of degrees thirteen and
  fourteen.
\newblock {\em Quart.\ J.\ Pure Appl.\ Math.} {\bf 29}, 224--249 (1897).

\bibitem{MilnorStasheff1974}
J.~W. Milnor and J.~D. Stasheff.
\newblock {\em Characteristic Classes}. Annals of Mathematics Studies~{\bf 76}.
\newblock Princeton University Press, Princeton, NJ, 1974.

\bibitem{Moebius1886}
A.~F. M\"obius.
\newblock Mittheilungen aus {M}\"obius' {N}achlass: {I}.\ {Z}ur {T}heorie der
  {P}oly\"eder und der {E}lementarverwandtschaft.
\newblock {\em Gesammelte Werke II} (F.~Klein, ed.),  515--559. Verlag von
  S.~Hirzel, Leipzig, 1886.

\bibitem{Royle1987}
G.~F. Royle.
\newblock The transitive groups of degree twelve.
\newblock {\em J.\ Symb.\ Comput.} {\bf 4}, 255--268 (1987).

\bibitem{Sarkaria1990}
K.~S. Sarkaria.
\newblock A generalized {K}neser conjecture.
\newblock {\em J.\ Comb.\ Theory, Ser.\ B\/} {\bf 49}, 236--240 (1990).

\bibitem{Sarkaria1991}
K.~S. Sarkaria.
\newblock A generalized van {K}ampen-{F}lores theorem.
\newblock {\em Proc.\ Am.\ Math.\ Soc.} {\bf 111}, 559--565 (1991).

\bibitem{SchulteWills1986}
E.~Schulte and J.~M. Wills.
\newblock Geometric realizations for {D}yck's regular map on a surface of genus
  $3$.
\newblock {\em Discrete Comput.\ Geom.} {\bf 1}, 141--153 (1986).

\bibitem{Sparla1997}
E.~Sparla.
\newblock {\em Geometrische und kombinatorische Eigenschaften
  trian\-gu\-lier\-ter Mannig\-faltigkeiten. {\rm Dissertation}}.
\newblock Shaker Verlag, Aachen, 1997, 132 pages.

\bibitem{Sparla1998}
E.~Sparla.
\newblock An upper and a lower bound theorem for combinatorial $4$-manifolds.
\newblock {\em Discrete Comput.\ Geom.} {\bf 19}, 575--593 (1998).

\bibitem{VolodinKuznetsovFomenko1974}
I.~A. Volodin, V.~E. Kuznetsov, and A.~T. Fomenko.
\newblock The problem of discriminating algorithmically the standard
  three-dimensional sphere.
\newblock {\em Russ.\ Math.\ Surveys\/} {\bf 29}, No. 5, 71--172 (1974).

\bibitem{Walkup1970}
D.~W. Walkup.
\newblock The lower bound conjecture for $3$- and $4$-manifolds.
\newblock {\em Acta Math.} {\bf 125}, 75--107 (1970).

\bibitem{Wang1995}
Z.~Wang.
\newblock Classification of closed nonorientable $4$-manifolds with infinite
  cyclic fundamental group.
\newblock {\em Math.\ Res.\ Lett.} {\bf 2}, 339--344 (1995).

\bibitem{Wilson1976}
S.~E. Wilson.
\newblock \emph{{N}ew {T}echniques for the {C}onstruction of {R}egular {M}aps.
  {\rm {D}issertation}}.
\newblock University of Washington, 1976, 194 pages.

\bibitem{Zubr1988}
A.~V. {\v{Z}ubr}.
\newblock Classification of simply-connected topological $6$-manifolds.
\newblock {\em Topology and Geometry -- Rohlin Seminar} (O.~Ya. Viro, ed.).
  Lecture Notes in Mathematics~{\bf 1346},  325--339. Springer-Verlag, Berlin,
  1988.

\end{thebibliography}

\bigskip
\bigskip
\medskip

\noindent
Ekkehard G.\ K\"ohler\\
Technische Universit\"at Berlin\\
Fakult\"at II - Mathematik und Naturwissenschaften\\
Institut f\"ur Mathematik, Sekr. MA 6-1\\
Stra\ss e des 17.\ Juni 136\\
D-10623 Berlin\\
{\tt ekoehler@math.tu-berlin.de}\\[5mm]

\noindent
Frank H.\ Lutz\\
Technische Universit\"at Berlin\\
Fakult\"at II - Mathematik und Naturwissenschaften\\
Institut f\"ur Mathematik, Sekr. MA 6-2\\
Stra\ss e des 17.\ Juni 136\\
D-10623 Berlin\\
{\tt lutz@math.tu-berlin.de}

\end{document}